\documentclass[12pt, a4paper]{article}

\usepackage[numbers,square,sort&compress]{natbib}

\usepackage{color} 
\usepackage{amsmath}
\usepackage{amsthm}
\usepackage{mathtools}
\usepackage {booktabs}
\usepackage{graphicx}

\newtheorem{theorem}{Theorem}

\newtheorem{proposition}{Proposition}
\newtheorem{lemma}{Lemma}

\theoremstyle{definition}

\newtheorem{remark}{Remark}

\newcommand{\CRHat}[1]{{ #1}^{\text{\scriptsize nc}}}

\begin{document}

\title{Guaranteed lower eigenvalue bounds for Steklov operators using conforming finite element methods}
\date{}
\author{ \normalsize{Taiga Nakano}
	    \thanks{Graduate School of Science and Technology, Niigata University,8050 Ikarashi 2-no-cho, Nishi-ku, Niigata City, Niigata, 950-2181, Japan e-mail: t-nakano@m.sc.niigata-u.ac.jp, xfliu@math.sc.niigata-u.ac.jp}
	   \and \normalsize{Qin Li} 
	   \thanks{School of Mathematics and Statistics, Beijing Technology and Business University, Beijing, 100048, P. R. China, e-mail: liqin@lsec.cc.ac.cn, yuemeiling@lsec.cc.ac.cn}
	   \and \normalsize{Meiling Yue$^{~\dagger}$}
	   \and \normalsize{Xuefeng Liu$^{~*}$ \footnote{Corresponding author} }
	   }
%

\maketitle

\abstract{For the eigenvalue problem of the Steklov differential operator, an algorithm based on the conforming finite element method (FEM) is proposed to provide guaranteed lower bounds for the eigenvalues.The proposed lower eigenvalue bounds utilize the {\em a priori} error estimation for FEM solutions to non-homogeneous Neumann boundary value problems,
which is obtained by constructing the hypercircle for the corresponding FEM spaces and boundary conditions.
Numerical examples demonstrate the efficiency of our proposed method.
}



\section{Introduction}\label{sec1}

To evaluate bounds of the eigenvalues for differential operators is a fundamental problem in numerical analysis. There are many approaches proposed to deal with the eigenvalue bounds, for example, the eigenvalue perturbation theories of Kato, the Lehamnn--Maehly--Goerisch theorem, the intermediate method, the homotopy method, etc.; Refer to \cite{boffi2010finite,sun2016finite} for surveys of FEM approaches to eigenvalue problems and \cite[Chapter 10]{nakao2019numerical_book} for a survey of methods with the purposes of explicit eigenvalue bounds.

For the study of lower eigenvalue bounds using the finite element method (FEM), there are two new approaches in the past decade.
\begin{enumerate}
\item [(1)]
The asymptotic analysis of lower eigenvalue bounds  tells that,   for many nonconforming FEMs, the approximate eigenvalues tend to the exact eigenvalues from below if the mesh is fine enough;
see, e.g., \cite{YZL-2010,lin2010stokes,HHL-2014,zhang2021asymptotic} and the references therein. However, since it is difficult to validate that the mesh size is small enough or not, one cannot obtain the lower eigenvalue bounds through the asymptotic lower bounds directly.
\item [(2)]
Another approach 
aims to provide {\em explicit bounds} for the eigenvalues.
Early results of Carstensen et al. \cite{CGa,CGe} and other groups \cite{SeVe, XXYY} require some {\em a priori} information of eigenvalues, for example,
the separation condition
or rough eigenvalue bounds for certain eigenvalues.
%
Fully computable explicit eigenvalue bounds without any additional conditions are proposed under  different approaches. In \cite{Liu2011,Liu,LO,YXLiu}, Liu utilizes the projection error-based technique to obtain explicit eigenvalue bounds. Such ideas trace back to the work of Birkhoff \cite{birkhoff}, Kikuchi \cite{kikuchi-liu} and Kobayashi \cite{Kobayashi2011,kobayashi-2015},  and are applied to  solving various eigenvalue problems \cite{XXL,hu2016guaranteed,liao2019optimal,liu2018explicit,zhang2021guaranteed,gallistl2021computational}.
{Canc\`es} et al. utilizes the residue error-based technique to bound the eigenvalues and the eigenfunctions \cite{CanDusMadStaVoh2017,CanDusMadStaVoh2018,Martin2020}. Recently, 
Carstensen et al. proposes new FEM schemes to provide direct lower eigenvalue bounds \cite{carstensen2020skeletal,carstensen2021guaranteed,carstensen2021direct}.
\end{enumerate}

\medskip

The Steklov eigenvalue problem is one of the important eigenvalue problems for differential operators; see \cite{BS,BRS,JLTL} for a  systematic introduction of background and applications.
Below is a short review of the numerical approaches to the eigenvalues of Steklov eigenvalue problems.
The qualitative error estimation by conforming FEM for Steklov eigenvalue problems are discussed in \cite{BO,CCMM,LLZ,LST}, based on which \cite{BZY,LYang,Xie} study more efficient algorithms such as two-grid and multilevel methods to solve Steklov eigenvalue problems. The {\em a posteriori}  error estimates with conforming FEM and nonconforming Crouzeix-Raviart FEM are discussed in \cite{AP} and \cite{DA}, respectively.
Especially, in \cite{LLX,YLL}, the asymptotic lower bounds for Steklov eigenvalue problems are discussed along with nonconforming finite elements. In \cite{YXLiu}, explicit lower bounds for the Steklov eigenvalues are obtained by using the Crouzeix-Raviart finite element along with an extension of the lower bound theorem of \cite{Liu}. 

\medskip

This paper considers the explicit eigenvalue bounds for the Steklov differential operator by using conforming FEMs. Since the positive semi-definite bilinear form $b(\cdot, \cdot)$ appears in the eigenvalue problem formulation $a(u,v)=\lambda b(u,v)$, we follow the theorem proposed in \cite{YXLiu} to handle the kernel space induced by $b(\cdot, \cdot)$.  
In \cite{YXLiu}, the following lower eigenvalue bound is proposed:
\begin{eqnarray*}
\lambda_{k}\geq \frac{\lambda_{k, h}}{1+M^{2}_{h}\lambda_{k, h}}.
\end{eqnarray*}
Here, $\lambda_{k, h}$'s are the non-infinity approximate eigenvalues obtained through the Galerkin method and $M_h$ is the quantity for the Galerkin projection error estimation.  
Different from the approach of \cite{YXLiu}, which utilized the nonconforming Crouzeix--Raviart FEM, this paper estimates $\lambda_{k,h}$ and $M_h$ by using the $H^1$-conforming FEM along with an explicit {\em a priori } error estimation for the projection operator.
Since the Neumann boundary condition is involved in the boundary value problem and the {\em a priori} error estimation has to deal with the worst case of the solution regularity, the obtained estimation of $M_h$ has the convergence  rate as $O(h^{1/2})$; see detailed discussion in Remark \ref{remark:discussion_convergence_rate}. 
The efficiency of proposed lower eigenvalue bounds is compared with the one of \cite{YXLiu} through numerical results. 
Since the conforming FEM is used in solving the eigenvalue problem, the approximate eigenvalue $\lambda_{k,h}$ gives upper bound of $\lambda_{k}$ directly.

\medskip

The rest of the paper is organized as follows. In section 2, we discuss the basic knowledge of the objective eigenvalue problem, its conforming linear finite element approximation and the lower eigenvalue bound theorem. 
In section 3, we show the details of the utilization of Hypercircle method in obtaining the explicit {\em a priori} error estimation, e.g., the value of $M_h$. In section 4, numerical results are shown to verify the theorem's results. We draw a conclusion in the last section.

\section{Objective eigenvalue problem and lower eigenvalue bounds}
Let $\Omega \subset \mathcal{R}^{2}$ be a bounded polygonal domain.
Throughout this paper, we use the standard notation (see, e.g. \cite{BaO,boffi2010finite}) for the Sobolev spaces $H^m(\Omega)$ ($m>0$).
Denote by $\|v\|_{L^2}$ or $\|v\|_{0}$ the $L^2$ norm of $v\in L^2(\Omega)$;
$|v|_{m,\Omega}$ and $\|v\|_{m,\Omega}$ the seminorm and norm in $H^{m}(\Omega)$, respectively.
Symbol $(\cdot, \cdot)$ denotes the inner product in $L^2(\Omega)$ or $(L^2(\Omega))^2$.
The space $H(\mbox{div}, \Omega)$ is defined by
$$
H(\mbox{div}, \Omega):=\{q\in (L^2(\Omega))^2\mid \mbox{div}\,q\in L^2(\Omega) \}.
$$

We are concerned with the following model Steklov eigenvalue problem:
\begin{equation}\label{eq:steklov_eig_pro_strong}
-\Delta u+ cu = 0\ \ \  {\rm in}\  \Omega\:;\quad
\frac{\partial u}{\partial \mathbf{n}}=\lambda u \ \ \ {\rm on}\ \Gamma=\partial\Omega \:,
\end{equation}
where $\frac{\partial }{\partial  \mathbf{n}}$ is the outward normal derivative on boundary $\Gamma$; $c$ is a non-negative number.

For a positive  $c$, we take $V=H^1(\Omega)$. In case $c=0$, the eigenvalue problem \eqref{eq:steklov_eig_pro_strong} has the zero eigenvalue and the eigenfunctions associated to the non-zero eigenvalues have  zero integral on the boundary of the domain. 
Upon this property of the eigenfunctions, let us 
take $V:=\{v\in H^1(\Omega) ~:~ \int_{\Gamma} v \mbox{d}s =0\}$ when $c=0$.

A weak formulation of the above problem is as follows: Find $\lambda\in\mathcal{R}$ and $u \in
V$ such that $\|u\|_{b}=1$ and
\begin{equation}
\label{eq:eig-pro-formulation}
a(u,v)= \lambda b(u,v) \quad \forall v\in V\:,
\end{equation}
where
$$
a(u,v):=\int_{\Omega}\nabla u\nabla v+cuv~\mbox{d}x,\quad
b(u,v):=\int_{\partial\Omega}uv~\mbox{d}s, \quad
\|u\|_{b}=\sqrt{b(u,u)}\:.
$$
Evidently the bilinear form $a(\cdot,\cdot)$ is symmetric,
continuous and coercive over $V$. 
The norm induced by $a(\cdot, \cdot)$ is denoted by $\|u\|_{a} := \sqrt{a(u,u)}$.

Let us consider the operator
$\mathcal{D}^{-1}:L^2(\Gamma) \rightarrow V$
such that for $f\in L^2(\Gamma)$, $\mathcal{D}^{-1}f=u$
satisfies the variational equation
\begin{eqnarray*}
a(\mathcal{D}^{-1}f, v)=b(f,v)\quad \forall v\in V.
\end{eqnarray*}
As a  compatibility  condition for the definition of $\mathcal{D}^{-1}$, it is required that  $\int_\Gamma f =0$ in case $c=0$. 
Let $\gamma$ be the trace operator $\gamma :V \rightarrow L^2(\Gamma)$.
Under the current assumption that the domain has a polygonal boundary, $\mathcal{D}^{-1}\circ\gamma :V \to V$ is a compact operator \cite{Demengel-2012}. 
The operator $\mathcal{D}^{-1}\circ\gamma$ has the zero eigenvalue, for which the associated eigenspace is just  $H^1_0(\Omega)$. The rest eigenvalues of $\mathcal{D}^{-1}\circ\gamma$ form a sequence $\{\mu_k\}$ as follows:
$$
\mu_k >0,\quad \mu_{1} \ge \mu_{2} \ge \cdots ,\quad  \lim\limits_{k\rightarrow \infty}\mu_{k}=0.
$$
In the rest of the paper, the trace operator $\gamma$ will be omitted if there is no ambiguity. 
The weak formulation of the eigenvalue problem for $\mathcal{D}^{-1}\circ\gamma$ is given by:  Find $u \in V$ and $\mu \ge 0$ such that,
\begin{equation}
\label{eq:new-formulation-eig}
b(u, v) = \mu a(u,v)\quad \forall v \in V\:.
\end{equation}
The eigenfunctions of  (\ref{eq:new-formulation-eig}) form a complete orthonormal basis of $V$.

As for the relation between the eigenvalue problem of $\mathcal{D}^{-1}\circ\gamma $ and the one defined in (\ref{eq:eig-pro-formulation}), we have that the {non-zero} eigenvalues $\mu_k$'s are given by the reverse of $\lambda_k$, i.e., $\mu_k={1}/{\lambda_k}$. 
\medskip

From the above argument, the eigenvalue problem
(\ref{eq:eig-pro-formulation}) has an eigenvalue sequence $\{\lambda_k \}:$
$$0<\lambda_1\leq \lambda_2\leq\cdots\leq\lambda_k\leq\cdots,\ \ \
\lim_{k\rightarrow\infty}\lambda_k=\infty~.
$$

\medskip

\noindent{\bf Finite element approximations}
Let $\mathcal{T}_{h}$ be a shape regular triangulation of the domain $\Omega$.
For each element $K\in \mathcal{T}_{h}$, denote by $h_{K}$ the longest edge length of $K$ and define the mesh size $h$ by the maximal value of $h_K$.
Particularly, it is assumed that, at  corners of the domain, 
each boundary edge of the triangulation is only shared by one triangle. Such an assumption is utilized in the proof of Lemma \ref{lem:solution_perturbation} to have a sharper error estimation.

The piecewise linear $H^1$-conforming finite element space $V^{h}$ is defined by
\begin{eqnarray*}
V^h&:=&\{v_{h}\in V:~v_{h}|_{K}\in P_1(K)~~ \forall K \in  \mathcal T_h \}~,
\end{eqnarray*}
where $P_1(K)$ is the space of polynomials of degree $\leq$ 1 on $K$.

The conforming finite element approximation of (\ref{eq:eig-pro-formulation}) is defined as follows: Find $\lambda_{h} (>0) \in\mathcal{R}$ and $u_{h} \in V^{h}$ such that 
$\|u_h\|_b=1$ and 
\begin{eqnarray}
\label{eq:eig-pro-FEMformulation}
a(u_{h}, v_{h})&=&\lambda_{h}b(u_{h},v_{h}) \quad\forall v_{h}\in V^{h}.
\end{eqnarray}
Let $n:=\mbox{dim}(V^{h})$ and $n_0:=n-\mbox{dim}(V^h\cap H_0^1(\Omega))$.
The eigenvalue problem (\ref{eq:eig-pro-FEMformulation}) has $n_0$ positive eigenvalues
$$
0<\lambda_{1,h}\leq\lambda_{2,h}\leq \cdots\leq\lambda_{n_0,h}<
 \infty \quad (n_0 \le n)~.
$$

\vskip 0.5cm


\medskip

Define the projection $P_h: V\rightarrow V^h$ by
$$
a(u-P_hu,v_h)=0 \quad \forall v_h\in V^h~.
$$

Below is the result from \cite{YXLiu} that provides lower eigenvalue bounds.
\begin{theorem}\label{2333}
Suppose the following inequality holds for the projection error:
\begin{equation*}
    \label{eq:Mh-def}
    \|(I-P_h) u\|_b \le M_h \|(I-P_h)u\|_a\quad \forall u \in V~.
\end{equation*}
Let $\lambda_{k,h}$ be the $k$-th  eigenvalue of \eqref{eq:eig-pro-FEMformulation}. 
A lower bound of the eigenvalue $\lambda_k$ of  \eqref{eq:eig-pro-formulation} is given by
\begin{equation}
\label{eq:final_eig_estimation}
\lambda_{k}\geq \frac{\lambda_{k, h}}{1+M^{2}_{h}\lambda_{k, h}},\quad k=1,\cdots, n_0.
\end{equation}
\end{theorem}

The algorithm to determine  the quantity $M_h$ with an explicit value is provided in the next section.

\section{Finite element approximation of the Neumann boundary value problem}\label{sec3}
The following boundary value problem and its FEM approach will play an important role in bounding the eigenvalues of the Steklov operator.
\begin{equation*}
\label{eq:neumann_bcv}
-\Delta u+ cu = 0 \mbox{ in }  \Omega; \quad 
\frac{\partial u}{\partial \mathbf{n}}=f \mbox{ on } \Gamma=\partial\Omega \:.
\end{equation*}
{Note that in case $c=0$, $f$ is further required to satisfy $\int_{\partial \Omega} f\mbox{d}s=0$.}

The weak formulation of the above problem is to find $u \in
V$ such that
\begin{equation}\label{eq:weak-form-bvp}
a(u,v)= b(f,v) \quad \forall v\in V.
\end{equation}
The conforming finite element approximation of (\ref{eq:weak-form-bvp}) is defined as follows: Find $u_{h} \in V^{h}$ such that
\begin{equation}\label{12}
a(u_{h}, v_{h})=b(f,v_{h}) \quad \forall v_{h}\in V^{h}.
\end{equation}

In this section, the following classical finite element spaces will be used in constructing the {\em a priori} error estimate for the FEM solution.
Let $E_{h}$ be the set of edges of the triangulation, and
$E_{h,\Gamma}$ the set of edges on the boundary of the domain.
Let  $\mathcal{T}^b_h$ be the set of elements of $\mathcal{T}_h$ having at least one edge on $\partial\Omega$.

\begin{itemize}
    \item [(i)]
Piecewise function spaces $X^{h}$ and $X_{\Gamma}^{h}$:
 \begin{eqnarray*}
 X^{h}&:=&\{v\in L^2(\Omega): v|_K\in P_1(K)\quad \forall K \in  \mathcal T_h\}\\
 X_{\Gamma}^{h}&:=&\{v\in L^2(\Gamma): v|_e\in P_1(e)\quad \forall e\in E_{h,\Gamma}\}
\end{eqnarray*}
where $P_1(e)$ is the space of polynomials of degree $\leq$ 1 on the edge $e$. {In case that $c = 0$, we further assume that $\int_{\Gamma} v \:\text{d}s =0$ for $v\in X_\Gamma^h$.}

\item [(ii)] The Raviart--Thomas FEM space $W^{h}$ with order one (\cite{BF}):
\begin{eqnarray*}
W^{h}:=\Big\{p_{h}\in H(\mbox{div}, \Omega)&\mid& p_{h}=
\left(
  \begin{array}{c}
    a_{K}\\
   b_{K}\\
  \end{array}
\right)+c_{K}
\left(
  \begin{array}{c}
    x\\
   y\\
  \end{array}
\right),\\
&&a_{K}, b_{K}, c_{K}\in P_{1}(K) ~
\mbox{for}~ K\in \mathcal{T}_{h} \Big\}.
\end{eqnarray*}
{ The freedoms of the Raviart-Thomas FEM space can be defined by the normal trace of $p_h$ on the edges of the triangulation. Hence, $\{ (p_h\cdot \mathbf{n}) |_{\Gamma}~|~p_h \in W^h\}=X^h$.  }
The space $W^{h}_{f_h}$ is a subset of $W^h$ corresponding to $f_{h}\in X^{h}_{\Gamma}$:
$$
W^{h}_{f_h}:=\{ p_{h}\in W^{h} \: | \: p_{h}\cdot \mathbf{n}=f_{h}~\mbox{on}~\Gamma \}.
$$
In {particular}, $W_0^h:=\{p_{h}\in W^{h} \: | \: p_{h}\cdot \mathbf{n}=0~\mbox{on}~\Gamma \}.$
\end{itemize}

Under current space settings, the following relations are available. 

\begin{equation*}
V^h \subset X^h, \quad \mbox{div} (W^h)=X^h, \quad \gamma  (V^h) \subset X_\Gamma^h~ .    
\end{equation*}

\subsection{The hypercircle method}

In this subsection, we introduce the hypercircle to be used to facilitate the error estimate in solving the eigenvalue problem. Let us introduce the following semi-norm (or norm if $c>0$) for $p \in H(\mbox{div};\Omega)$:
$$
\|p\|^2_{H(\text{div}),c}:=\int_{\Omega} |\mbox{div}\; p|^2 + c |p|^2 \mbox{d}\Omega~.
$$

\begin{theorem}
\label{thm:hypercircle}
Given $f_h \in X_\Gamma^h$, 
let $u$ be the solution of (\ref{eq:weak-form-bvp}) with $f:=f_h$. For $v_h\in V^h$ and $p_{h}\in W_{f_h}^h$ satisfying $\text{div } p_h =c v_h$,  the following hypercircle holds:
\begin{eqnarray}\label{eq:hypercircle}
\|u-v_h\|_{a}^2+
\|\nabla u-p_h\|_{H(\text{div}),c}^2=\|\nabla v_h-p_h\|_{L^2}^2~.
\end{eqnarray}
\end{theorem}
\proof
Rewriting $\nabla v_h-p_h$ by $(\nabla v_h-\nabla u)+(\nabla u-p_h)$, we have
\begin{eqnarray*}
\|\nabla v_h-p_h\|_{L^2}^2=\|\nabla v_h-\nabla u\|_{L^2}^2+
\|\nabla u- p_h\|_{L^2}^2+2(\nabla v_h-\nabla u,\nabla u-p_h).
\end{eqnarray*}
Furthermore, the Green theorem and the Neumann boundary conditions setting lead to
\begin{eqnarray*}
(\nabla u_h-\nabla u,\nabla u-p_h) &=&  ( v_h-u, -c u+\mbox{div } p_h)\\
&=& ( v_h-u,-cu+cv_h) =c\|u-v_h\|_{L^2}^2.
\end{eqnarray*}
Noticing that $\|\nabla u-p_h\|_{H(\text{div}),c}^2=\|\nabla u-p_h\|_{L^2}^2+c\|u- v_h\|_{L^2}^2$, we obtain the hypercircle in \eqref{eq:hypercircle}.
\endproof

Next, let us introduce the quantity $\kappa_{h}$ such that
\begin{eqnarray}
\label{eq:def-kappa-h}
\kappa_{h}&:=&\max_{f_{h}\in X_{\Gamma}^{h}\setminus \{0\}} ~
\min_{\substack{v_h \in V^h,~ p_h \in W_{f_h}^h \\ \text{div}\: p_h =cv_{h} }} ~
\frac{\|\nabla v_h-p_h\|_{L^2}}{\|f_{h}\|_{b}}.
\end{eqnarray}

\begin{lemma}
\label{lem:error_est_using_kappa_h}
Given $f_h \in X_\Gamma^{h}$, let $\tilde{u}\in
V$ and $\tilde{u}_{h} \in V^{h}$ be the solutions to the following variational problems, respectively,
\begin{eqnarray}
&& a(\tilde{u},v)=b(f_{h},v) \quad\forall v\in V, \label{21} \\
&&
\nonumber a(\tilde{u}_{h},v_{h})= b(f_{h},v_{h}) \quad\forall v_{h}\in V^{h}. 
\end{eqnarray}
Then, the following error estimate holds:
\begin{eqnarray}\label{eq:computable estimation}
\|\tilde{u}-\tilde{u}_{h}\|_{a}\leq \kappa_{h}\|f_{h}\|_{b}~.
\end{eqnarray}
\end{lemma}
\proof
In Theorem \ref{thm:hypercircle}, take $v_{h}:= \tilde{u}_{h}$,  $u:=\tilde{u}$ and $p_h\in W_{f_h}^h$ such that $\mbox{div\:}p_h =c\tilde{u}_{h}$, then we have
\begin{eqnarray}\label{eq:a-posteriori-est}
\|\tilde{u}-\tilde{u}_{h}\|_{a}\leq \|\nabla\tilde{u}_{h}-p_{h}\|_{L^2}.
\end{eqnarray}
By further considering the minimization of $p_h$ and the variation of $f_h$ in $X_{\Gamma}^h$, we draw the conclusion in (\ref{eq:computable estimation}).
\endproof

\begin{remark}
In Theorem 3.3 of \cite{LLiu}, a general case such that  $\mbox{div }p_h - c \tilde{u}_h \not=0$ is discussed, for which the formulation of $\kappa_h$ is little complicated with a free parameter to be adjusted properly. Since the Raviart--Thomas space $W^h$ in this paper has a higher order, one can find $p_h\in W^h$ such that $\mbox{div }p_h =  c\tilde{u}_h$ holds for {$\tilde{u}_h$} $\in V^h$. As a defect of the current setting, the  Raviart--Thomas space $W^h$ with a higher order will cause larger matrices in the computation. 
In \eqref{eq:new_kappa_h} of \S \ref{sec:kappa_h_computing}, 
a new quantity $\bar{\kappa}_h$, which can be solved with improved computation efficiency, is proposed to produce a reasonable upper bound of $\kappa_h$.
\end{remark}

\subsection{Explicit {\em a priori} error estimates}

We first quote an explicit bound for the constant in trace theorem.
A direct estimation of $C_{e}(K)$ with FEM approximations is also provided in \S \ref{sec:direct_estimate_c_k}.
\begin{lemma}[\cite{YXLiu}]\label{lemma:perturbation_of_f}
Let $e$ be an edge of triangle element $K$.
Define function space
$$
V_{e}(K):=\{v\in H^{1}(K) ~~|~~ \int_{e}v ~\text{d}s=0 \}\:.
$$
Given $u\in V_e(K)$, we have the following inequality related to the trace theorem:
\begin{equation}
\label{eq:c_lliu}
  \|u\|_{L^{2}(e)}\leq   C_{e}(K)|u|_{H^{1}(K)},
\quad    C_{e}(K) := 0.574 \sqrt{ \frac{|e|}{|K|}}h_{K} \le 0.8118 \frac{h_{K}}{\sqrt{H_K}} \:.
\end{equation}
Here, $H_K$ denotes the height of triangle $K$ with respect to edge $e$.
\end{lemma}

Given an element $K$ of $\mathcal{T}^h$ with $e$ as one of its edges,  
let $\pi_{0,e}$ be the linear operator that takes the average of a function on edge $e$. Let $I$ be the identity operator. Note that $\pi_{0,e} v$ is defined over the element $K$.
For function $v\in H^1(\Omega)$,  $(I-\pi_{0,e})v|_K$ is regarded as a shift of $v$, that is,
$$
(I-\pi_{0,e}) v|_K = v|_K - \frac{1}{|e|}\int_{e} v \mbox{d}s \in H^1(K) ~.
$$
Since $(I-\pi_{0,e}) v|_K$ has zero integral on the boundary edge $e$, the following error estimation holds:
\begin{equation}
    \label{eq:pi_0_error_est}
\|(I-\pi_{0,e}) v\|_{L^2(e)} \le C_{e}(K) |v|_{H^1(K)}~.
\end{equation}

{
Let us introduce a piecewise $L^2$ projection operator $\pi_{h,\Gamma}: L^2(\Gamma)\mapsto X_{\Gamma}^{h}$ on the boundary faces}: Given $f\in L^2(\Gamma)$,
$\pi_{h,\Gamma}f \in X_{\Gamma}^{h}$ satisfies
$$
b(f-\pi_{h,\Gamma}f, v_{h})=0 \quad \forall v_{h}\in X_{\Gamma}^{h}.
$$
It is easy to see that on a boundary edge $e$ of $\mathcal{T}^h$,
$$
\int_e (f-\pi_{h,\Gamma} f)|_e \; \pi_{0,e} v\; \mbox{d} s=0\quad \forall v \in H^1(\Omega)~.
$$

\medskip

\begin{lemma}\label{lem:solution_perturbation}
Let $u$ and $\tilde{u}$
be solutions to (\ref{eq:weak-form-bvp}) and (\ref{21}), respectively, with $f_h$ taken as $f_h:=\pi_{h,\Gamma}f$. Then, the following error estimate holds:
\begin{eqnarray}
\label{eq:sol_pert_respect_to_f}
\|u-\tilde{u}\|_{a}
\leq C_{e,h}\|(I-\pi_{h,\Gamma})f\|_{b},
\end{eqnarray}
where $C_{e,h}$ takes the maximum of $C_{e}(K)$ over the boundary elements:
\begin{equation*}
C_{e,h}:=\max \limits_{K \in \mathcal{T}^b_h} C_{e}(K) = O(h^{1/2})\:.
\end{equation*}
\end{lemma}
\proof
Setting $v=u-\tilde{u}$ in (\ref{eq:weak-form-bvp}) and (\ref{21}), we have
\begin{eqnarray}
&& \nonumber a(u-\tilde{u}, u-\tilde{u}) \\
 \notag\\
&=&b(f-f_{h},u-\tilde{u}) =\sum_{e\subset E_{h, \Gamma}} \int_e (I-\pi_{h,\Gamma})f \cdot (I-\pi_{0,e})(u-\tilde{u}) \mbox{d} s \notag \\
&\le & \|(I-\pi_{h,\Gamma})f\|_{b}
\left\{\sum_{e\in E_{h, \Gamma}} \|(I-\pi_{0,e})(u-\tilde{u})\|_{L^2(e)}^2\right\}^{1/2} \label{eq:local-trace-est-1}
.
\end{eqnarray}
By applying the estimation \eqref{eq:pi_0_error_est}, we have 
\begin{equation}
\label{eq:local-trace-est-2}
\sum_{e\in E_{h, \Gamma}} \|(I-\pi_{0,e})(u-\tilde{u})\|_{L^2(e)}^2 
\le
\sum_{K \in \mathcal{T}^b_h}
C_{e}(K)^2 |u-\tilde{u}|_{H^{1}(K)}^2
\le 
C_{e,h}^2 \|u-\tilde{u}\|_{a}^2~.
\end{equation}
Note that, the first inequality of the above estimation holds under the assumption that each boundary edge of the triangulation is only shared by one triangle.
For a general mesh without such an assumption, the coefficient in the estimation should be doubled.
The  estimations \eqref{eq:local-trace-est-1} and \eqref{eq:local-trace-est-2} lead to the estimation (\ref{eq:sol_pert_respect_to_f}). The convergence rate of $C_{e,h}$ as $C_{e,h}=O(h^{1/2})$ for regular meshes is obvious from the estimation \eqref{eq:c_lliu}.

\endproof

Now, we are ready to propose the explicit {\em a priori} error estimation.

\begin{theorem}
\label{thm:a-priori-est}
Let $u$ and $u_{h}$
be solutions to (\ref{eq:weak-form-bvp}) and (\ref{12}), respectively. The following error estimates hold.
\begin{equation}
\label{eq:a_priori_error_estimation}
\|u-u_{h}\|_{a} \leq  M_{h}\|f\|_{b},~~ \|u-u_{h}\|_{b}\leq M_h \|u-u_{h}\|_{a}  \leq  M_{h}^2\|f\|_{b},
\end{equation}
where  $M_{h}:=\sqrt{C_{e,h}^2+\kappa_{h}^2}.$
\end{theorem}
\proof
Take $f_h := \pi_{h,\Gamma}f$ and consider the decomposition $f=f_h +(f-f_h)$. Let $\tilde{u}_h$ be the one defined in Lemma \ref{lem:solution_perturbation} corresponding to $f_h$. The minimization principle for the FEM solution $u_h$ tells that
$\|u-u_{h}\|_{a}\leq \|u-\tilde{u}_{h}\|_{a}$. 
By further applying (\ref{eq:computable estimation}) of Lemma \ref{lem:error_est_using_kappa_h}  
 and \eqref{eq:sol_pert_respect_to_f} of Lemma \ref{lem:solution_perturbation}, we have
\begin{eqnarray*}
\|u-u_{h}\|_{a}&\leq& \|u-\tilde{u}_{h}\|_{a}\leq  \|u-\tilde{u}\|_{a}+\|\tilde{u}-\tilde{u}_{h}\|_{a}\\
&\leq& C_{e,h}\|(I-\pi_{h,\Gamma})f\|_{b}+\kappa_{h}\|f_{h}\|_{b}\\
&\leq& \sqrt{C_{e,h}^2+\kappa_{h}^2}\|f\|_{b}=M_h\|f\|_{b}\:.
\end{eqnarray*}

The error estimate (\ref{eq:a_priori_error_estimation}) can be obtained by applying the standard Aubin--Nitsche duality technique.
\endproof

\begin{remark}
\label{remark:discussion_convergence_rate}
{ The analysis of $C_{e,h}$ tells that $C_{e,h}=O(h^{1/2})$,  and numerical results in \S\ref{subsec:kappa_computation} imply that $\kappa_h$ has the convergence rate as $O(h^{1/2})$  even for convex domains and high-order FEM spaces. }
Hence, the proposed {\em a priori} error estimation with the quantity $M_h$ has the convergence rate as $O(h^{1/2})$, which will lead to a lower eigenvalue bound given by \eqref{eq:final_eig_estimation} with a degenerated  convergence rate as $O(h)$.
From classical discussions of the solution regularity of Neumann boundary condition,
%
{ it is known that the solution has the regularity as $u\in H^{1+r}(\Omega)$ for a general $f\in L^2(\partial\Omega)$ with $r\in [0,1/2)$; see, e.g., \cite[Theorem 4]{savare1998regularity} 
and \cite[Theorem 31.34]{guermond2021finite}.} 
Therefore, such a convergence rate of $M_h$ is reasonable, as the {\em a priori} error estimation has to manipulate the worst case of the  solution regularity. 
Meanwhile, the FEM approximations of the leading eigenvalues over the unit square domain demonstrate the $O(h^2)$ convergence rate  (see the discussion in Section 4). 
{
Thus, as the defect of the proposed lower eigenvalue bounds in this paper, the estimation \eqref{eq:final_eig_estimation} using  $M_h=O(h^{1/2})$ is sub-optimal for smooth eigenfunctions. }
\end{remark}
{
\begin{remark}
It is worth pointing out that Theorem \ref{thm:a-priori-est} is also available for general $\mathcal{R}^n$  $(n\ge 2)$ spaces by providing explicit values for the involved quantities. The value of $\kappa_h$ can be computed by using the hypercircle for standard FEM spaces on $\mathcal{R}^n$ domain. 
For the constant $C_{e,h}$ appearing in Lemma \ref{lemma:perturbation_of_f}, the method used in \cite{YXLiu} to evaluate $C_{e,h}$ can be easily extended to a $\mathcal{R}^n$ simplex; see such a discussion in, e.g., the corrigendum of \cite[Lemma 1]{AINSWORTH2014184}.
\end{remark}
}

\subsection{Computation of $\kappa_{h}$}
\label{sec:kappa_h_computing}

This section is dedicated to a description of the algorithm to evaluate $\kappa_{h}$ defined in (\ref{eq:def-kappa-h}). 

\medskip
First, for a fixed $f_h \in X_\Gamma^h$, 
we consider the following minimization problem:
\begin{equation*}
    \min_{u_h \in V^h} \min_{\substack{p_h \in W_{f_h}^h \\ \text{div}\: p_h =c u_{h}}} {\|\nabla u_h-p_h \|_{L^2}^2}~.
\end{equation*}
The above problem is reformulated as finding the stationary point for the following objective function: for $(u_h, p_h,x_h)  \in V^h \times W^h_{f_h} \times X^h$, 
\begin{equation*}
\mathcal{F}(u_h, p_h,x_h) :=
  \frac{\|\nabla u_h-p_h \|_{L^2}^2}{2} + (x_h,\text{div}\:p_h-c u_h). 
\end{equation*}
Then, stationary point $({u}_{h}, p_h,x_h)$ satisfies
\begin{equation}
\label{eq:v1}
    \left\{
    \begin{array}{lllll}
        & (\nabla u_h, \nabla v_h) & -(p_h, \nabla v_h)  & -c(x_h, v_h) &= 0 \\
        &- (\nabla u_h, q_h)&+ (p_h, q_h) &+(x_h, \text{div}\:q_h) &= 0\\
        &- c(u_{h}, y_h)&+(\text{div}\:p_h, y_h) &\quad &= 0
    \end{array}
    \right.
\end{equation}
for all $(v_h,q_h,y_h) \in V^h \times W^h_{0} \times X^h$. 

To confirm the existence and uniqueness of $(u_h, p_h,x_h)$ of the system \eqref{eq:v1}, we cite the following result from \cite{BF}. {Note that the notation below is restricted to the discussion of Proposition \ref{thm:suddle point} in the rest of current subsection.}

\begin{proposition}[Proposition 1.1 of \cite{BF},  p.38]
\label{thm:suddle point}
Let $V$ and $Q$ be Hilbert spaces, the dual spaces of which are denoted by $V'$ and $Q'$, respectively. Let $B:V \to Q'$ be an linear operator.
Let $g \in \mathrm{Im}(B)$ and let the bilinear form $a(\cdot,\cdot)$ be coercive on $\mathrm{Ker}(B)$, that is, there exists $\alpha_0$ such that 
 $$
  a(v_0, v_0) \ge \alpha_0 \|v_0\|^2 \quad \forall v_0 \in \mathrm{Ker}(B).
 $$
 Then, given $f\in V'$, there exists a unique $u \in V$ solution of the equations:
 $$
  Bu = g;\quad
  a(u, v_0) = \langle f,v_0 \rangle _{V'\times V }\quad \forall v_0 \in \mathrm{Ker}(B)
  ~.
 $$
\end{proposition}


To apply Proposition \ref{thm:suddle point}, we
consider a reformulation of  (\ref{eq:v1}).
Let  $\hat{p}_{h}$ be a fixed function of $W^h_{f_h}$ and introduce $p_{h,0}:=p_h - \hat{p}_{h} \in W^h_0$.  
The equations in  (\ref{eq:v1}) becomes
\begin{equation*}
\label{eq:v2}
    \left\{
    \begin{array}{lllll}
        & (\nabla u_h, \nabla v_h) &-(p_{h,0}, \nabla v_h)  & - c(x_h, v_h) &= (\hat{p}_{h}, \nabla v_h) \\
        &- (\nabla u_h, q_h)&+ (p_{h,0}, q_h) &+(x_h, \text{div}\:q_h) &= -(\hat{p}_{h}, q_h)\\
        & - c(u_h, y_h)&+(\text{div}\:p_{h,0}, y_h) &\quad &= -(\text{div}\: \hat{p}_{h}, y_h)
    \end{array}
    \right.~.
\end{equation*}
Let us consider the following function settings.
\begin{align*}
&V := V^h \times W_0^h,\quad Q := X^h,\\
&\langle f, \{v_h, q_h\} \rangle_{V' \times V} := (\hat{p}_{h}, \nabla v_h-q_h)_\Omega,\quad 
\langle g, \cdot \rangle_{Q' \times Q} := (-\text{div}\: \hat{p}_{h},\cdot)_\Omega~,
\\ &a(\{u_h,p_{h,0}\},\{v_h,q_h\})
:= (\nabla u_h -p_{h,0}, \nabla v_h-q_h)_\Omega, 
\\
&\langle  B(\{u_h,p_{h,0}\}), \cdot \rangle_{Q' \times Q}  := (\text{div}\: p_{h,0} -cu _h,  \cdot)_\Omega~.
\end{align*}
The inner product of $V$ is defined by
$$
\langle \{u_h, p_{h}\}, \{v_h, q_h\} \rangle_V
:=
(\nabla u_h, \nabla v_h) + c (u_h, v_h) + (p_{h}, q_h) + (\text{div } p_{h}, \text{div } q_h) ~,
$$
which induces the norm as 
$\|\{u_h,p_h\}\|_{V} =\{\|\nabla u_h\|_{\Omega}^2 + c \|u_h\|_{\Omega}^2 + \|p_h\|_{H(\text{div})}^2 \}^{\frac{1}{2}}$.
Since the involved spaces are finite dimensional,  $\mathrm{Im}(B)$ is the closed subspace of  $V^h \times W_0^h$.
The positive-definiteness and boundedness of $a(\cdot,\cdot)$ are easy to confirm.

The coercivity of $a(\cdot, \cdot)$
over $\mathrm{Ker}(B)$ can be confirmed by the following equality:
for $\{u_h, p_{h,0}\} \in \mathrm{Ker}(B)$,  by applying Green's formula,
\begin{eqnarray*}
 a(\{u_h,p_{h,0}\},\{u_h,p_{h,0}\}) 
 & = & \|\nabla u_h \|^2 - 2(\nabla u_h, p_{h,0}) + \|p_{h,0} \|^2 \\
 & = & \|\nabla u_h \|^2 + 2c\|u_h\|^2 + \|p_{h,0} \|^2 \\
 & = & \|\nabla u_h \|^2 + c\|u_h\|^2  + \|\text{div } p_{h,0}\|^2 + \| p_{h,0} \|^2 \\
 & = & \|\{u_h, p_{h,0}\}\|_{V}^2.
\end{eqnarray*}

\medskip

Therefore, Proposition \ref{thm:suddle point} makes certain that the functional $\mathcal{F}$ has a unique saddle point $(u_h, p_{h,0} + \hat{p}_h, x_h)$ in $V^h\times W_{f_h}^h\times X^h$, giving a solution to the problem. 
The evaluation of $\kappa_h$ can be done by further considering the maximization of ${\|\nabla u_{h}-p_h \|_{L^2}^2}/{\|f_h\|_{b}}^2$ for all $f_h \in X_\Gamma^h$.

\medskip

In the practical computation, we propose an efficient way that provides an upper bound for $\kappa_h$.
Given an $f_h\in X_\Gamma^h$, let us consider the following formulation that determines $\tilde{u}_{h}\in V^{h}$ and $p_{h}\in W^{h}_{f_h}$ subsequently. 
\begin{itemize}
\item [(a)] Find $\tilde{u}_{h}\in V^{h}$ s.t.
\begin{eqnarray*}\label{59}
a(\tilde{u}_{h}, v_{h})=b(f_{h},v_{h})\quad \forall v_{h}\in V^{h}.
\end{eqnarray*}
\item [(b)] Let $\tilde{u}_{h}$ be the solution of (a).
~Find $p_{h}\in W^{h}_{f_h}$ and $\rho_{h}\in X^{h}$, $r\in \mathcal{R}$ s.t.
\begin{equation*}
\left\{
\begin{array}{rcll}
(p_{h},{q}_{h})+(\rho_{h},\mbox{div } {q}_{h}) + (\rho_h, s) &=& { 0 } & \forall q_{h}\in W_0^{h}, ~ {\forall s \in \mathcal{R}} \\
(\mbox{div } p_{h},\eta_h) + ({r}, \eta_h) &=& c(\tilde{u}_{h}, \eta_h)\quad & \forall \eta_h\in X^{h}
\end{array}
\right.~.
\end{equation*}

\end{itemize}

For each given $f_h$, there exist unique solution $\tilde{u}_{h}$ and $p_{h}$ to the sub-problems (a) and (b).
By using the mapping from $f_h$ to $\tilde{u}_{h}$ and $p_{h}$,
let us introduce the quantity $\bar{\kappa}_{h}$,
which works as an upper bound of $\kappa_{h}$:
\begin{equation}
\label{eq:new_kappa_h}
\bar{\kappa}_{h}:=\max_{f_{h}\in X_{\Gamma}^{h}\setminus \{0\}} ~
\frac{\|\nabla\tilde{u}_{h}-p_{h}\|_{0}}{\|f_{h}\|_{b}}.
\end{equation}
According to the definition of $\bar{\kappa}_h$, it is required to find $f_{h}$ that maximizes the
value of $\|\nabla\tilde{u}_{h}-p_{h}\|_{0}/\|f_{h}\|_{b}$, which can be achieved by solving an eigenvalue problem for matrices. Since 
$\tilde{u}_{h}\in V^{h}$ and $p_{h}\in W^{h}_{f_h}$ are determined subsequently, the matrices involved in setting up the linear system will has a quite smaller size than the ones in solving \eqref{eq:v2}. 
For detailed description of the evaluation of $\kappa_h$ and $\tilde{\kappa}_h$, refer to (\cite{LO}), where an analogous problem is considered.

\begin{remark} The introduction of variable $r$ in the setting of problem (b) is to make certain a regular matrix in solving the linear systems. By setting $v_{h}=1$ in the problem (a), we have
$$c\int_{\Omega}\tilde{u}_{h}~\text{d}x=\int_{\partial\Omega}f_h ds =\int_{\partial\Omega}p_h\cdot\mathbf{n}~\text{d}s=
\int_{\Omega}\mbox{div } p_{h}~\text{d}x\:.$$
The above relation implies that 
$
(\mbox{div }p_h -c \tilde{u}_h, \cdot) 
$
has a kernel space with constant function.
\end{remark}

\section{Numerical Examples}\label{sec6}

In this section, we apply the eigenvalue estimation (\ref{eq:final_eig_estimation}) along with the explicit {\em a priori} error estimation solve the eigenvalue problem (\ref{eq:steklov_eig_pro_strong}) on both the unit square domain $\Omega=(0,1)\times (0,1)$ and the L-shaped domain
$\Omega=(0, 2)\times (0, 2)\setminus [1, 2]\times [1, 2]$. 
Here, we select $c$ appearing in (\ref{eq:steklov_eig_pro_strong}) as 1.
Also, the existing method of \cite{YXLiu} based on the nonconforming FEM is utilized to compare the efficiency with each other.

\subsection{Evaluation of $\kappa_h$  and $\bar{\kappa}_h$}
\label{subsec:kappa_computation}
We adopt two different methods in subsection \ref{sec:kappa_h_computing} to evaluate $\kappa_h$ and $\bar{\kappa}_h$ and display the computation results in Tab. \ref{table0}-\ref{table01}. It is observed that the
$\bar{\kappa}_h$ gives very close upper bound of $\kappa_h$;
for the square domain, 
the leading $4$ significant digits
of $\bar{\kappa}_h$ and $\kappa_h$ are the same to each other. Thus, $\bar{\kappa}_h$ will be utilized instead of $\kappa_h$ in the following computation examples.
It is worth to point out that the value of $\kappa_h$ has a convergence rate, denoted by $\gamma(\kappa_h)$ in the tables, as $O(h^{1/2})$ for both the square domain and the L-shaped domain. { To confirm the dependency of the convergence rate of $\kappa_h$ on the order of FEM spaces, the hypercircle using FEM spaces (i.e., $V^h, W^h, X^h, X^h_{\Gamma}$) of order 2 
is used to evaluate $\kappa_h$, denoted by $\kappa_{h,2}$, is also displayed in Table \ref{table0}. Numerical results tell that $\gamma(\kappa_{h,2})$ is still $0.5$.}

\begin{table}[ht]
\begin{center}
\caption{Quantities $\kappa_h, \bar{\kappa}_h$ and $\kappa_{h,2}$ for the unit square domain ($\gamma$: convergence rate)}\label{table0}
\begin{tabular}{ccccc}
$h$&$\sqrt{2}/4$&${\sqrt{2}}/{8}$&${\sqrt{2}}/{16}$&${\sqrt{2}}/{32}$ \\ \midrule
$\kappa_h$& 0.2891 & 0.2042  &0.1443 & 0.1021  \\
~ $\gamma(\kappa_h)$ ~ & - & 0.50 & 0.50 & 0.50  \\
\hline
$\bar{\kappa}_h$ &0.2891 & 0.2042&0.1443&0.1021\\
~ $\gamma(\bar{\kappa}_h)$ ~ & - & 0.50 & 0.50 & 0.50 \\
\hline
$\kappa_{h,2}$& 0.2291 & 0.1621  &0.1146 & 0.0811  \\
~ $\gamma(\kappa_{h,2})$ ~ & - & 0.50 & 0.50 & 0.50  \\
\end{tabular}
\end{center}
\end{table}
\begin{table}[ht]
\begin{center}
\caption{Quantities  $\kappa_h$ and $\bar{\kappa}_h$ for the L-shaped domain domain ($\gamma$: convergence rate) }\label{table01}
\begin{tabular}{ccccc} 
$h$&$\sqrt{2}/2$&${\sqrt{2}}/{4}$&${\sqrt{2}}/{8}$&${\sqrt{2}}/{16}$ \\ \midrule
$\kappa_h$& 0.5075 & 0.3624 & 0.2588 & 0.1846 \\
~ $\gamma({\kappa}_h)$ ~ & - & 0.49 & 0.49 & 0.49 \\ 
\hline
$\bar{\kappa}_h$&0.5106 & 0.3633 & 0.2591 & 0.1847\\
~ $\gamma(\bar{\kappa}_h)$ ~ & - & 0.49 & 0.49 & 0.49 \\ 
\end{tabular}
\end{center}
\end{table}

\begin{figure}
    \centering
    \includegraphics[width=3.3cm]{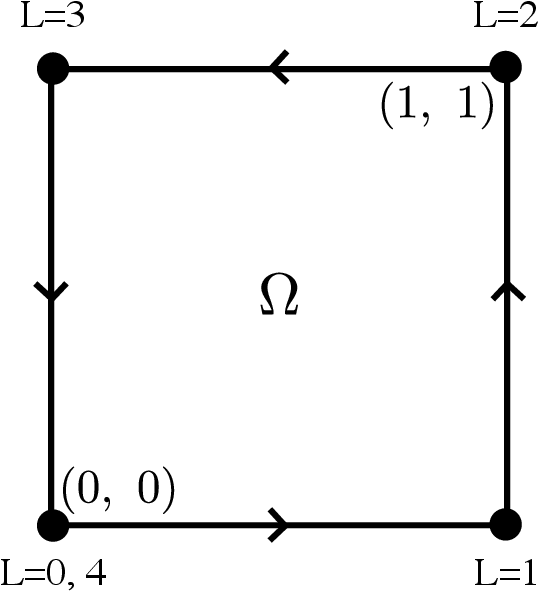} \quad
    \includegraphics[width=3.3cm]{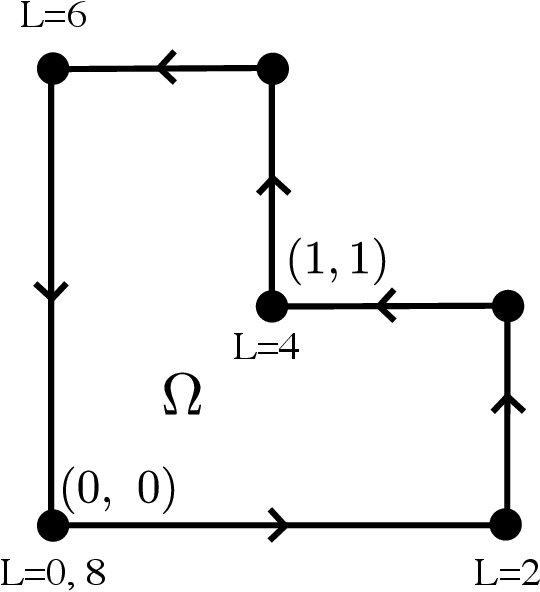}
    \caption{Parameter $L$ for the arc length of domain boundary}
\label{fig:arc_length_of_domain_boundary}
\end{figure}

\begin{figure}[ht]
    \centering
    \includegraphics[height=3.3cm]{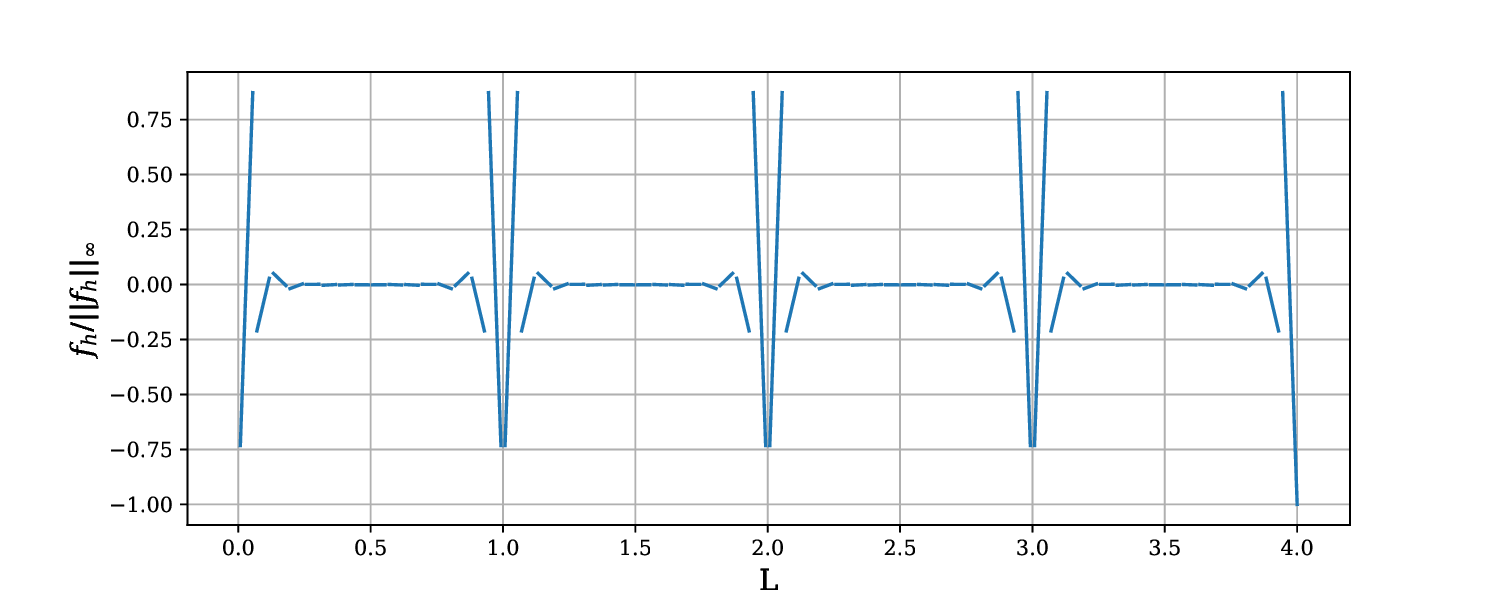}
    \includegraphics[bb=0 0 432 288, height=3.3cm]{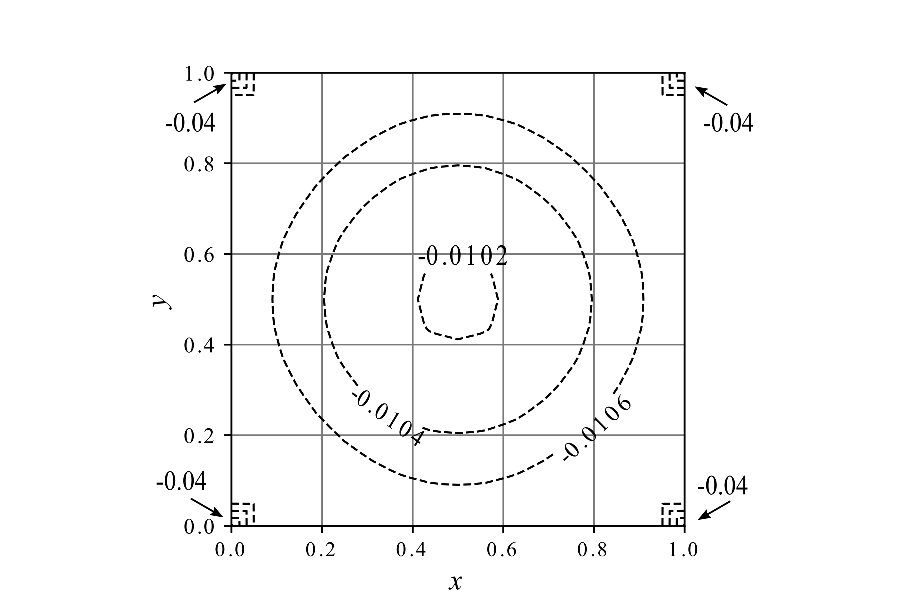}
    \caption{The worst $f_h$ (left) and $u_h$ (right) that determine $\kappa_h$ (square domain)}
\label{fig:worst_f_u_square}
\end{figure}

\begin{figure}[h!]
    \centering
    \includegraphics[height=3.3cm]{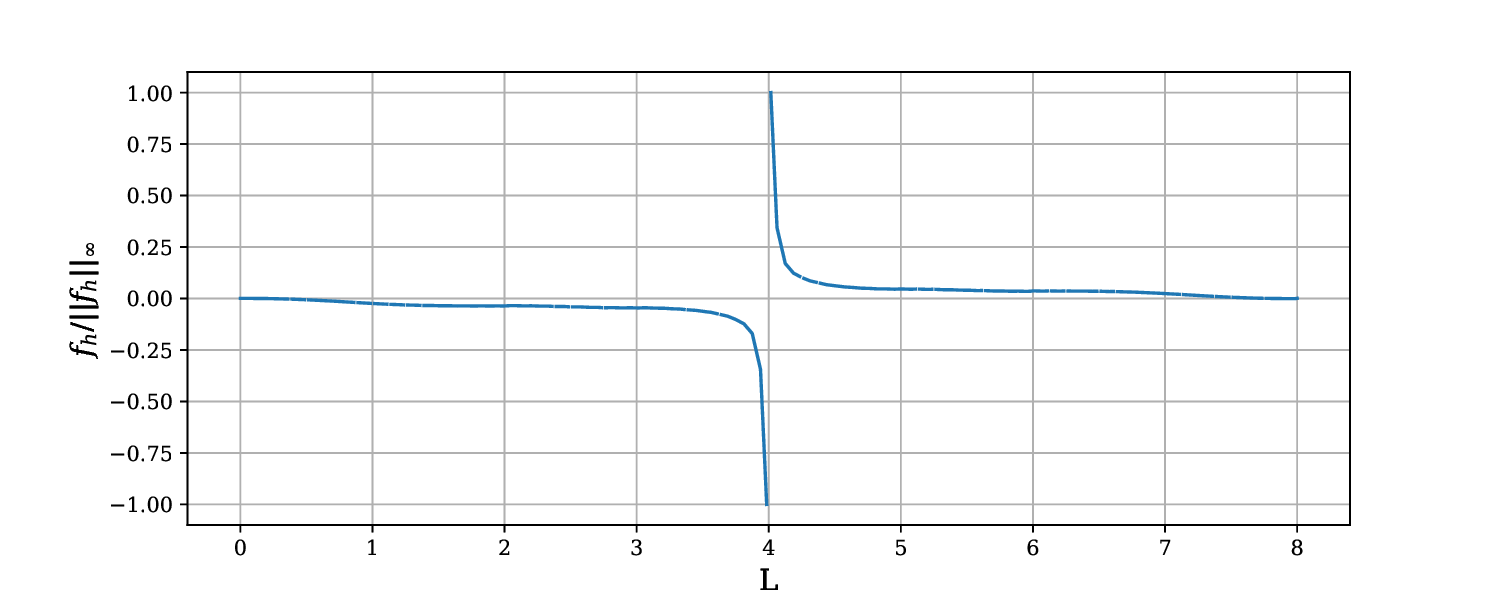} \includegraphics[height=3.3cm]{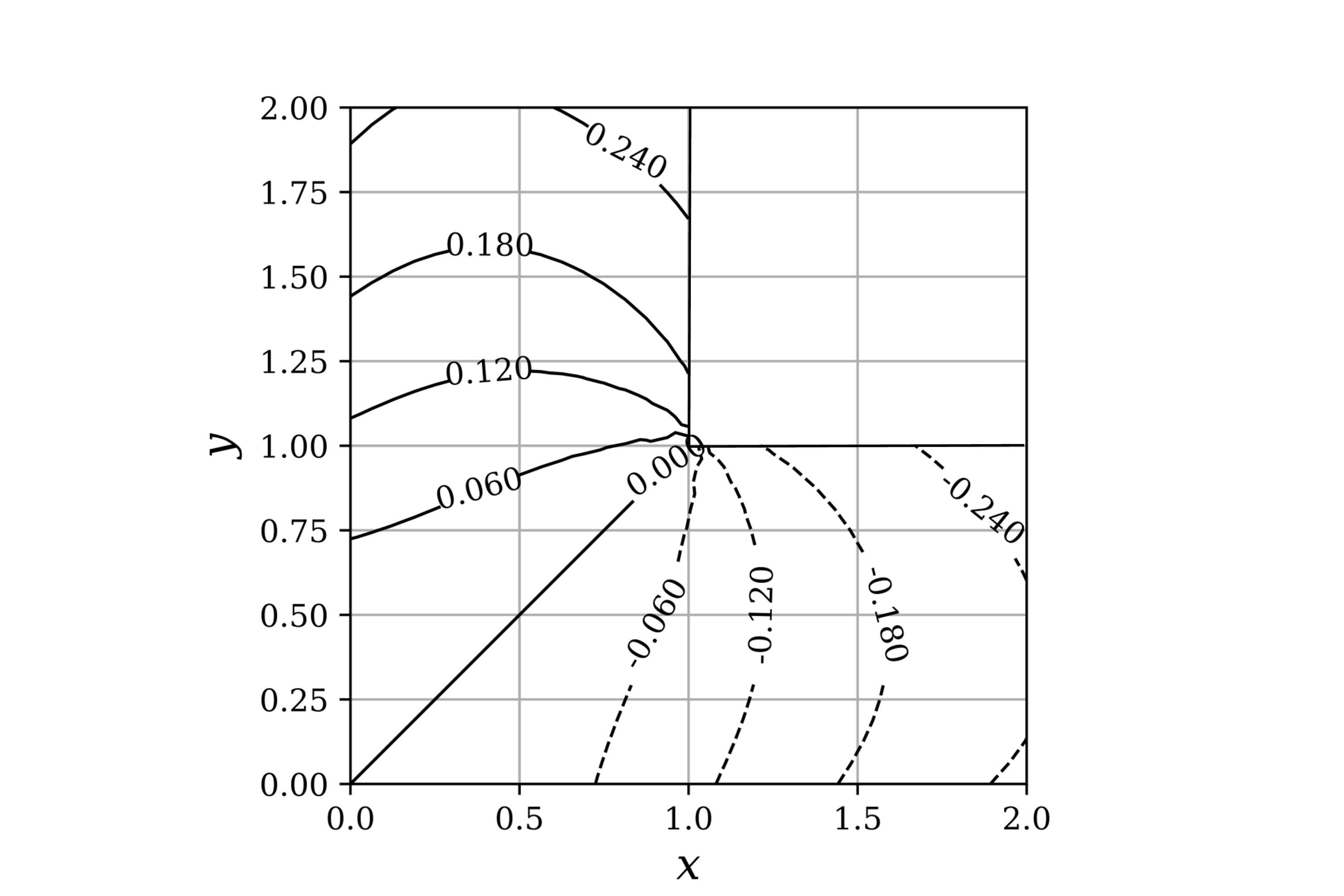}
    \caption{The worst $f_h$ (left) and $u_h$ (right) that determine $\kappa_h$ (L-shaped domain)}
\label{fig:worst_f_u_l_shape}\end{figure}

It is of great interest when the worst case of the projection error happens. To confirm for which 
$f_h$ the value of $\kappa_h$ is reached, we draw the figures of such an $f_h$ and its corresponding conforming FEM solution $u_h$.
Since $f_h$ is defined on the boundary of domain, let us introduce a parameter $L$ to measure the arc length from the vertex located at the origin point; see Fig. \ref{fig:arc_length_of_domain_boundary}. The graphs of $f_h$ and the contour lines of $u_h$ for the square domain and the L-shaped domain are displayed in Fig. \ref{fig:worst_f_u_square} and \ref{fig:worst_f_u_l_shape}, respectively. Note that
$f_h$ is normalized by the $L^\infty$ norm in each figure.
The numerical results imply that when the value of $f$ is concentrated at the corner of the domain, the worst case of the projection error happens. For the square domain, there is large variation of both $f_h$ and the conforming FEM solution $u_h$ around the four corners, while for the L-shaped domain, the variation of both $f_h$ and $u_h$ is concentrated at the re-entry corner. A theoretical investigation of the worst cases for the Neumann boundary conditions is of interest and will be considered in the future work.

\subsection{Preparation for eigenvalue estimation}

\begin{figure}[hb!]
\centering
\includegraphics[width=6cm]{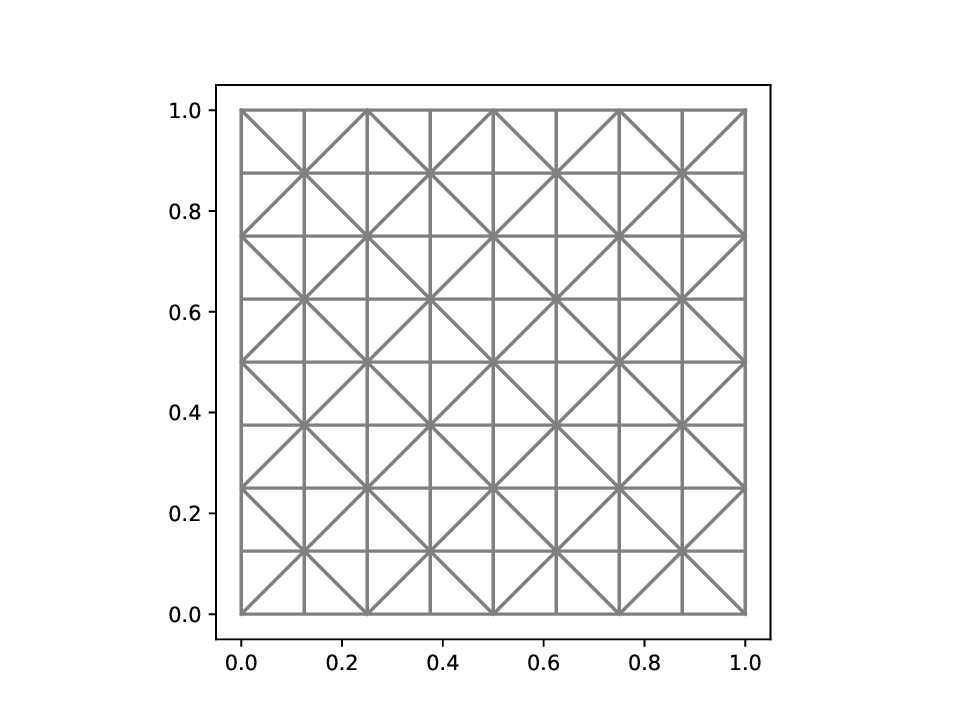}
\includegraphics[width=6cm]{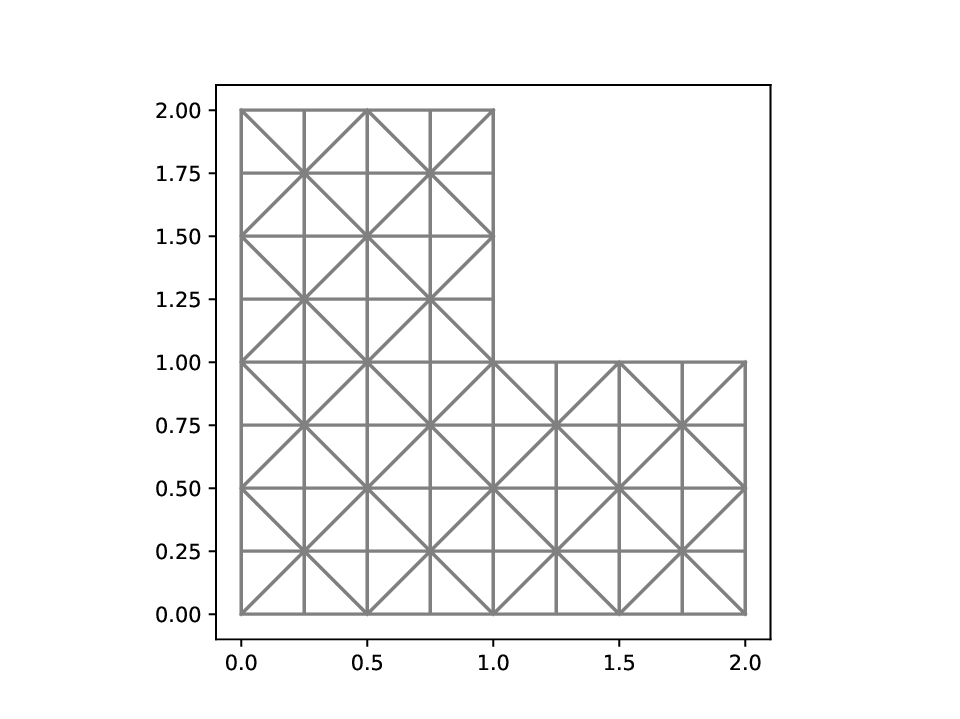}
\vspace{-0.45cm}
\caption{The unit square and L-shaped domains}
\label{fig:0}
\end{figure}

The explicit values of the exact eigenvalues for both domains are not available.
For the unit square domain, the following high-precision estimation with reliable significant digits are used as a nice approximation
to true eigenvalues (\cite{YLL}).
$$
(\mbox{unit square})~\lambda_1\approx 0.240079, \quad \lambda_2= \lambda_3 \approx 1.49230\:.
$$
In case of the L-shaped domain, the cubic conforming FEM with the mesh size $h=\sqrt{2}/256$ provides a high-precision approximation to eigenvalues:
$$
(\mbox{L-shaped domain})~ \lambda_1\approx 0.3414160,  \quad \lambda_2 \approx 0.6168667, \quad \lambda_3 \approx 0.9842784\:.
$$

For both domains, the uniform meshes are adopted.  The eigenvalue estimation (\ref{eq:final_eig_estimation}) provides a guaranteed lower eigenvalue bound:
\begin{equation}
\label{eq:CG-eigenbound}
\underline{\lambda}_{k, h} := \frac{\lambda_{k, h}}{1+M_{h}^{2}\lambda_{k, h}}, \quad \quad M_h=\sqrt{C_{e,h}^2+\kappa_{h}^2}~,
\end{equation}
where $\lambda_{k, h}$ denotes the $k$-th approximate eigenvalue from the conforming FEM and
the quantity $C_{e,h}$ in estimating $M_{h}$ is given by
$$
C_{e,h} := 0.8118\max\limits_{K\in \mathcal{T}^b_h}\frac{h_K}{\sqrt{H_K}} (=0.9654\sqrt{h_K})\:.
$$
{Note that $h_K = \sqrt{2} H_K$.}
The eigenvalue estimation from Theorem 3.8 of \cite{YXLiu} has the formula as follows.
\begin{equation}
\label{eq:eig_estimation_you_xie_liu}
\CRHat{\underline{\lambda}}_{k, h}:=\frac{\CRHat{\lambda}_{k, h}}{1+\widehat{C}^2_{h}\CRHat{\lambda}_{k, h}}\:,
\end{equation}
where $\CRHat{\lambda}_{k, h}$ denotes the $k$-th approximate eigenvalue from the Crouzeix-Raviart FEM.
Particularly, for the uniform mesh used here, $\widehat{C}_{e,h}$ is estimated by
\begin{eqnarray*}
\widehat{C}_{e,h}&= &0.6711\max\limits_{K\in \mathcal{T}^b_h}\frac{h_K}{\sqrt{H_K}}+\frac{0.1893}{\sqrt{\CRHat{\lambda}_{1,h}}}\max\limits_{K\in\mathcal{T}_h }h_K\\
&=&0.7981\sqrt{h_K}+\frac{0.1893}{\sqrt{\CRHat{\lambda}_{1,h}}}h_K~.
\end{eqnarray*}

\subsection{Computation results for two domains}
Sample uniform triangular meshes for two domains are displayed in Fig. \ref{fig:0}, where
the mesh size for the unit square is $h=\sqrt{2}/8$ and the one for the L-shaped domain is $h=\sqrt{2}/4$.

\begin{figure}[ht!]
\centering
\includegraphics[width=6.5cm]{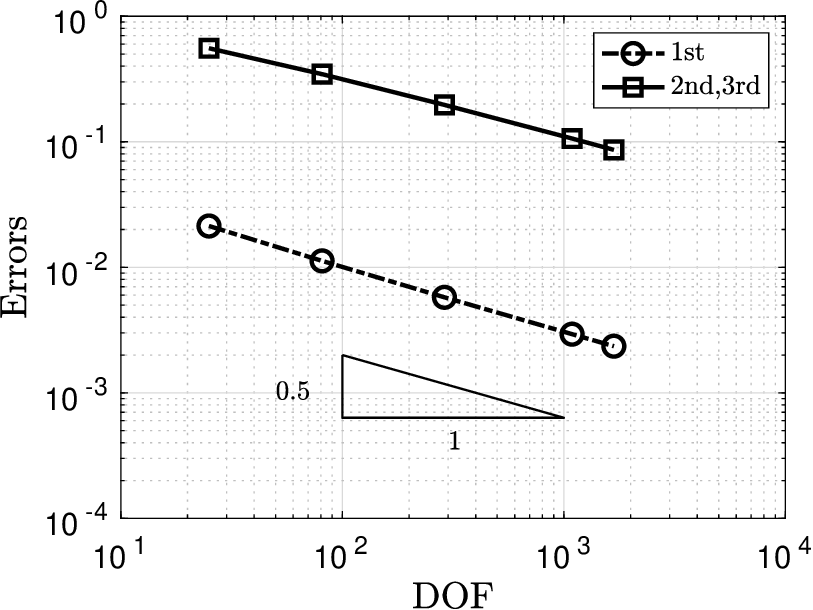} 
\includegraphics[width=6.5cm]{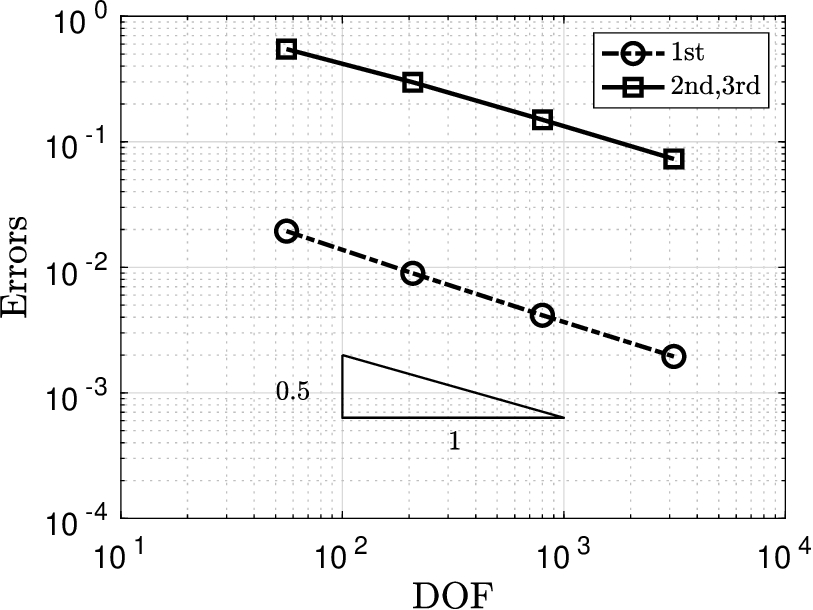}~~ 
\caption{Errors of eigenvalue bounds v.s. DOF (the unit square domain) 
(Left: $|\lambda_i-\underline{\lambda}_{i,h}|$ ; Right: $|\lambda_i-\CRHat{\underline{\lambda}}_{i,h}|$ ($i=1, 2, 3$))
}
\label{fig:1-2}
\end{figure}

\begin{figure}[h!]
\centering
\includegraphics[width=6.5cm]{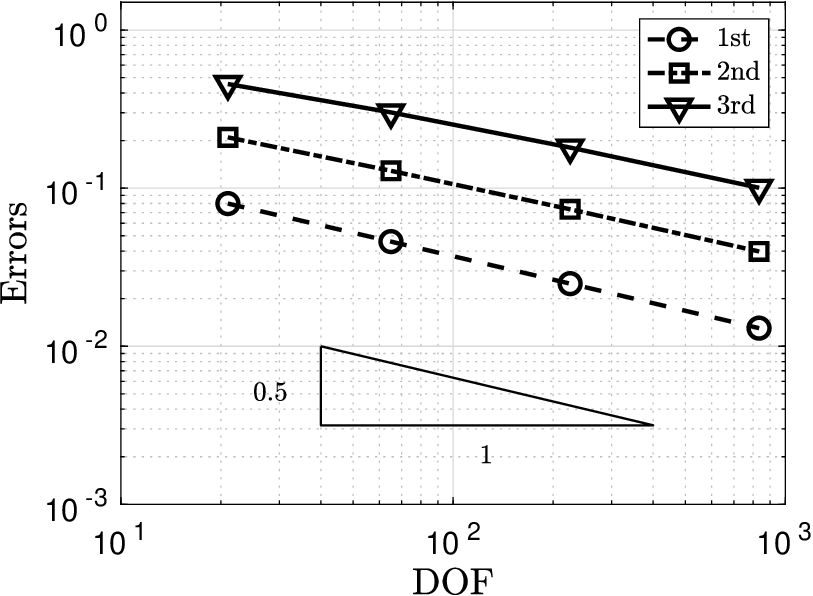} 
\includegraphics[width=6.5cm]{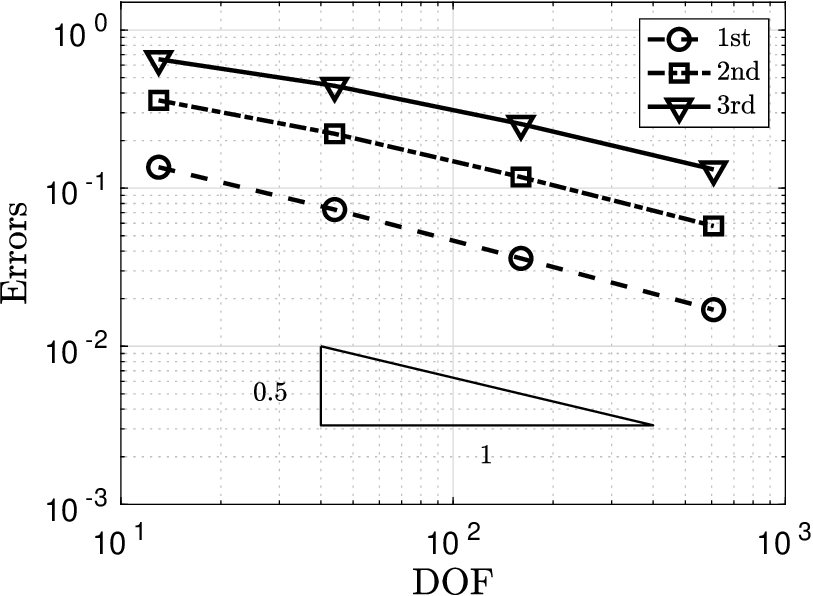}~~ 
\vspace{-0.3cm}
\caption{Errors of eigenvalue bounds v.s. DOF (the L-shaped domain) 
(Left: $|\lambda_i-\underline{\lambda}_{i,h}|$, Right: $|\lambda_i-\CRHat{\underline{\lambda}}_{i,h}|$ ($i=1, 2, 3$))
}
\label{fig:1-3}
\end{figure}

For the unit square domain, the eigenvalue estimations (\ref{eq:final_eig_estimation}) for the leading $3$ eigenvalues are displayed in Tab. \ref{table1},
while the results based on the nonconforming FEM (\cite{YXLiu}) are displayed in Tab. \ref{table2}.
The results for the L-shaped domain are displayed in Tab. \ref{table3} and \ref{table4}. Fig. \ref{fig:1-2} and Fig. \ref{fig:1-3}
describe the relation between the absolute errors and the degrees of freedom (DOF) over the unit square and  L-shaped domains, respectively. Here, the DOF of
(\ref{eq:final_eig_estimation}) is counted as the
the dimension of the linear conforming FEM space $V^h$,
while the one for \cite{YXLiu} is the dimension of the Crouzeix-Raviart FEM space.

\begin{figure}[h!]
\centering
\includegraphics[width=6.5cm]{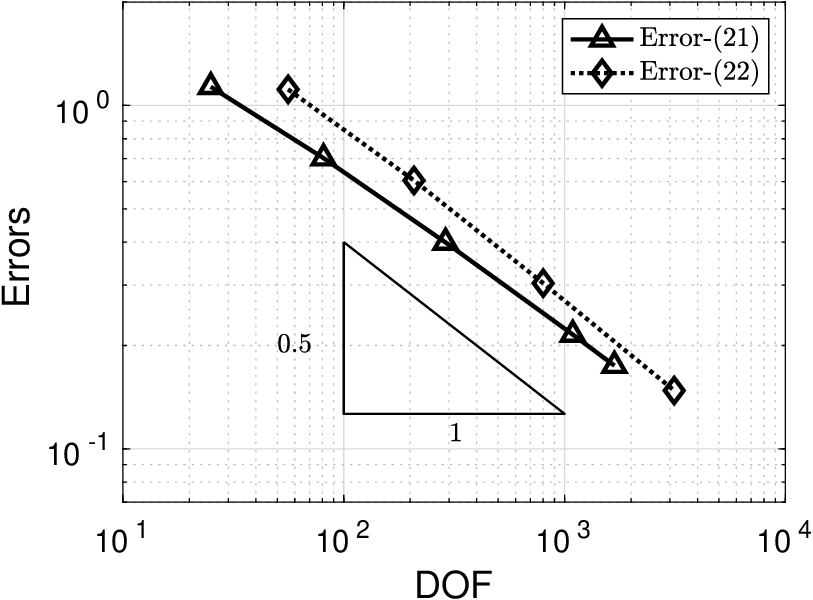}~~ 
\includegraphics[width=6.5cm]{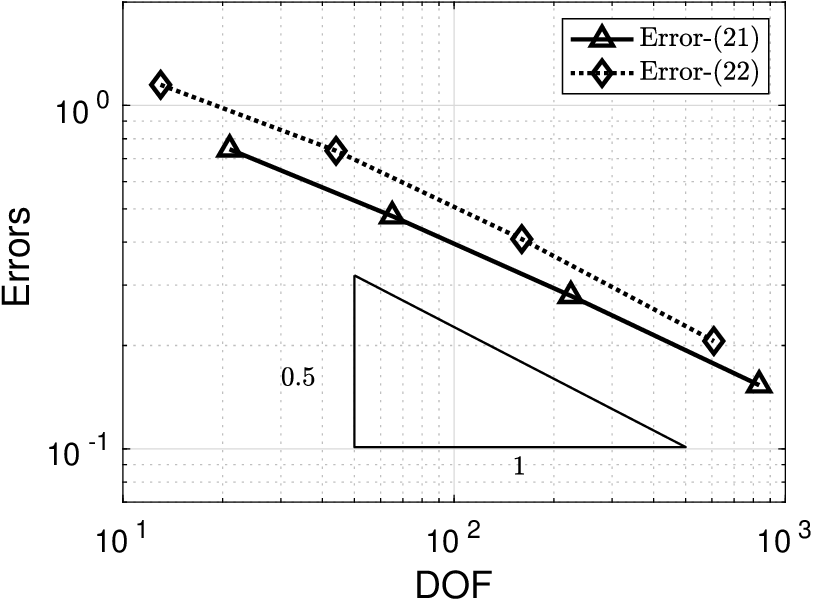} 
\caption{The total errors for the eigenvalue bounds v.s. DOF  (Left: the unit square; Right: the L-shaped domain) 
}
\label{fig:1}
\end{figure}

Let us also introduce the total errors by
\begin{eqnarray*}
\mbox{Error-}(\ref{eq:CG-eigenbound})&:=&|\lambda_1-\underline{\lambda}_{1,h}|+|\lambda_2-\underline{\lambda}_{2,h}|+|\lambda_3-\underline{\lambda}_{3,h}|~,\\
\mbox{Error-(\ref{eq:eig_estimation_you_xie_liu})}&:=&|\lambda_1-\CRHat{\underline{\lambda}}_{1,h}|+|\lambda_2-\CRHat{\underline{\lambda}}_{2,h}|+|\lambda_3-\CRHat{\underline{\lambda}}_{3,h}|~.
\end{eqnarray*}
The relation between the total errors and the degrees of freedom is displayed in Fig. \ref{fig:1}.

Different from the nonconforming FEM in \cite{YXLiu} which merely provide the guaranteed lower eigenvalue bounds, the conforming FEM produces both the upper bounds and the lower bounds of the eigenvalues.
From the computational results for two domains and the comparison between
the bound (\ref{eq:final_eig_estimation}) and the one from \cite{YXLiu}, we  draw the conclusion that
\begin{enumerate}
  \item [(1)]
  Both the lower eigenvalue bounds proposed in this paper and the one in \cite{YXLiu} have a sub-optimal convergence rate for the leading Steklov eigenvalues,
  compared with the convergence rate estimated by the numerical results themselves. 
  \item [(2)]
  With the same degree of freedom, the lower bound in (\ref{eq:final_eig_estimation}) (or \eqref{eq:CG-eigenbound})
  gives slightly better estimation than the one from the nonconforming FEM.
  However, to obtain the bound (\ref{eq:final_eig_estimation}), one has to pay more efforts to solve a matrix problem to obtain $\bar{\kappa}_h$. 
\end{enumerate}

\newpage

\begin{table}[ht!]
\begin{center}
\caption{Quantities in the eigenvalue estimation (\ref{eq:CG-eigenbound}) ($\gamma$: convergence rate; unit square domain)}\label{table1}
\begin{tabular}{c|cccc|c}
\rule[-2mm]{0mm}{6mm}{}$h$&$\sqrt{2}/4$&${\sqrt{2}}/{8}$&${\sqrt{2}}/{16}$&${\sqrt{2}}/{32}$ & $\gamma$ \\ \midrule 
\rule[-2mm]{0mm}{6mm}{}
$\bar{\kappa}_h$ &0.2891 & 0.2042&0.1443&0.1021&0.51 \\
\rule[-2mm]{0mm}{4mm}{}
$C_{e,h}$ &0.5740 & 0.4059&0.2870&0.2029&0.50  \\
\rule[-2mm]{0mm}{4mm}{}
$M_{h}$ &  0.6427 &   0.4544 &   0.3208 &   0.2272&0.51  \\
\hline
\rule[-2mm]{0mm}{6mm}{}
$\lambda_{1,h}$& 0.2404841  & 0.2401798  &0.2401042& 0.2400854& 2.01  \\
\rule[-2mm]{0mm}{4mm}{}
$\underline{\lambda}_{1, h}$ &  0.218753 & 0.228833  &0.2343144 &0.2371468 & 0.95  \\
\hline
\rule[-2mm]{0mm}{6mm}{}
$\lambda_{2,h}$& 1.527151 & 1.502305 &1.494918  & 1.492966 & 1.92 \\
\rule[-2mm]{0mm}{4mm}{}
$\underline{\lambda}_{2, h}$ &  0.936415 & 1.146662  & 1.295596& 1.386153 & 0.72\\
\end{tabular}
\end{center}
\hspace{2cm} (Note: $\lambda_{2,h}=\lambda_{3,h}$, $\underline{\lambda}_{2, h}=\underline{\lambda}_{3, h}$)
\end{table}

\begin{table}[ht!]
\begin{center}
\caption{Quantities in the eigenvalue estimation (\ref{eq:eig_estimation_you_xie_liu}) ($\gamma$: convergence rate; unit square domain)\label{table2} }
\begin{tabular}{c|cccc|c}
\rule[-2mm]{0mm}{6mm}{}$h$&${\sqrt{2}}/{4}$&${\sqrt{2}}/{8}$&${\sqrt{2}}/{16}$&${\sqrt{2}}/{32}$ & $\gamma$\\ \midrule 
\rule[-2mm]{0mm}{6mm}{}$\widehat{C}_{e,h}$&0.6110176 & 0.4038323  & 0.2714162 &0.1848489  &0.61\\  \hline
\rule[-2mm]{0mm}{6mm}{}
$\CRHat{\lambda}_{1,h}$&  0.2404829  &0.2401793  & 0.2401041    &  0.2400853&2.0
\\
\rule[-2mm]{0mm}{4mm}{}
$\CRHat{\underline{\lambda}}_{1, h}$& 0.2206705 & 0.2311264 &0.235931  &  0.2381318& 1.13 \\
\hline
\rule[-2mm]{0mm}{6mm}{}$\CRHat{\lambda}_{2,h}$&1.460229  &  1.483297&  1.489892  &  1.491678& 1.88
\\
\rule[-2mm]{0mm}{4mm}{}
$\CRHat{\underline{\lambda}}_{2, h}$&0.9450309  &1.19438  & 1.342541  & 1.419335&0.95 \\
\end{tabular} \\
\end{center}
\hspace{2cm} (Note: $\CRHat{\lambda}_{2,h}=\CRHat{\lambda}_{3,h}$, $\CRHat{\underline{\lambda}}_{2, h}=\CRHat{\underline{\lambda}}_{3, h}$)
\end{table}

\begin{table}[h]
\begin{center}
\caption{Quantities in the eigenvalue estimation {(\ref{eq:CG-eigenbound})} ($\gamma$: convergence rate; L-shaped domain) }\label{table3}
\begin{tabular}{c|cccc|c}
\rule[-2mm]{0mm}{6mm}{}$h$&${\sqrt{2}}/{2}$&${\sqrt{2}}/{4}$&${\sqrt{2}}/{8}$&${\sqrt{2}}/{16}$ & $\gamma$ \\ \midrule 
\rule[-2mm]{0mm}{6mm}{}$\bar{\kappa}_h$&0.5106&0.3633&0.2591&0.1847 & 0.48\\
\rule[-2mm]{0mm}{4mm}{}
$C_{e,h}$&0.8118&0.5740  & 0.4059 &0.2870 & 0.50  \\
\rule[-2mm]{0mm}{4mm}{}
$M_{h}$&0.9590 &0.6793  & 0.4815 & 0.3413 & 0.50 \\ \hline
\rule[-2mm]{0mm}{6mm}{}
$\lambda_{1,h}$&0.3443305& 0.3421498&0.3416010&0.3414626 &  2.06\\
\rule[-2mm]{0mm}{4mm}{}
$\underline{\lambda}_{1, h}$ &0.2615119&0.2954914 & 0.3165279  & 0.3283997 & 0.93 \\
\hline
\rule[-2mm]{0mm}{6mm}{}
$\lambda_{2,h}$ &0.6513041 &0.6299816&0.6217140&0.6186763 & 1.45 \\
\rule[-2mm]{0mm}{4mm}{}
$\underline{\lambda}_{2, h}$& 0.4073133&0.4880800  &  0.5433766 & 0.5770854 & 0.89\\
 \hline
\rule[-2mm]{0mm}{6mm}{}
$\lambda_{3,h}$&1.0278736 &0.9968693&0.9876317&0.9851393 & 2.02\\
\rule[-2mm]{0mm}{4mm}{}
$\underline{\lambda}_{3, h}$&  0.5283698  & 0.6827630  &  0.8035932  & 0.8837230 & 0.85 \\
\end{tabular}
\end{center}
\end{table}

\begin{table}[h!]
\begin{center}
\caption{Quantities in the eigenvalue estimation {(\ref{eq:eig_estimation_you_xie_liu})} ($\gamma$: convergence rate; L-shaped domain) \label{table4}}
\begin{tabular}{c|cccc|c} 
\rule[-2mm]{0mm}{6mm}{}
$h$&${\sqrt{2}}/{2}$&${\sqrt{2}}/{4}$&${\sqrt{2}}/{8}$&${\sqrt{2}}/{16}$ & $\gamma$\\ \midrule 
\rule[-2mm]{0mm}{6mm}{}$\widehat{C}_{e,h}$ &0.8997886 & 0.5890361  & 0.3928155&0.2659045&0.63     \\   \hline
\rule[-2mm]{0mm}{6mm}{}$\CRHat{\lambda}_{1,h}$ & 0.3425959  & 0.3416846  &0.3414799&  0.3414316& 2.08 \\
\rule[-2mm]{0mm}{4mm}{}$\CRHat{\underline{\lambda}}_{1, h}$  &0.2682036 &0.3054704&0.3243874 &0.3333834&1.07 \\
\hline
\rule[-2mm]{0mm}{6mm}{}$\CRHat{\lambda}_{2,h}$  &0.5829704  & 0.6039094 & 0.6120116& 0.6150436& 1.42   \\
\rule[-2mm]{0mm}{4mm}{}$\CRHat{\underline{\lambda}}_{2, h}$ & 0.3960439  &0.4992908  &0.5592028 &0.5894119&  0.99   \\
 \hline
\rule[-2mm]{0mm}{6mm}{}$\CRHat{\lambda}_{3,h}$ &0.9608929  &0.9769290  &0.9821661& 0.9837098&  1.76 \\
\rule[-2mm]{0mm}{4mm}{}$\CRHat{\underline{\lambda}}_{3, h}$  & 0.5404476&0.7296185 &0.8529063 &0.9197389& 0.88 \\
\end{tabular}
\end{center}
\end{table}

\medskip

\subsection{Comparison with the optimal $C_{e}(K)$ and proposed bound in (\ref{eq:c_lliu})}
\label{sec:direct_estimate_c_k}

In this subsection, we estimate 
the trace constant $C_{e}(K)$ over several triangle $K$'s directly, and compare with its bound in (\ref{eq:c_lliu}). 
For $i=1,2,3$, denote the $i$th edge of $K$ by $e_i$.
Let us introduce the function space $V_{e_i}$ after $V_{e}$ in Lemma \ref{lemma:perturbation_of_f}.
$$
V_{e_i}(K) = \{ u \in H^1(K) | \int_{e_i} u \:ds = 0 \}.
$$
The trace constant $C_{e_i}(K)$ is the quantity that makes certain the following estimation holds.
$$\|u\|_{L^{2}(e_i)}\leq   C_{e_i}(K)|u|_{H^{1}(K)} \quad \forall u \in V_{e_i}(K)~.
$$ 
The determination of $C_{e_i}(K)$ reduces to finding the minimal positive eigenvalue of the following Steklov eigenvalue problem:
\begin{equation}
\label{traceSteklov}
-\Delta u =0 \mbox{ in } K,~~ \frac{\partial u}{\partial \mathbf{n}} =\lambda u \mbox{ on } e_i, ~~~ \frac{\partial u}{\partial \mathbf{n}} = 0
\mbox{ on } \partial K \setminus e_i.
\end{equation}
By taking $a(u,v) := (\nabla u, \nabla v)_{K}, ~ b(u,v) := (u,v)_{e_i}$, the weak formulation of \eqref{traceSteklov} is given as follows:
$$
\mbox{Find } (\lambda, u) \in \mathcal{R} \times V_{e_i} \mbox{ s.t. }~~ a( u, v ) = \lambda b(u,v)\quad \forall v \in V_{e_i}(K).
$$
The strict lower eigenvalue bound for the above eigenvalue problem can be obtained by an analogous argument as performed in this paper, the detail of which is omitted here.

\begin{figure}[ht!]
\begin{center}
    \includegraphics[width=3.6cm]{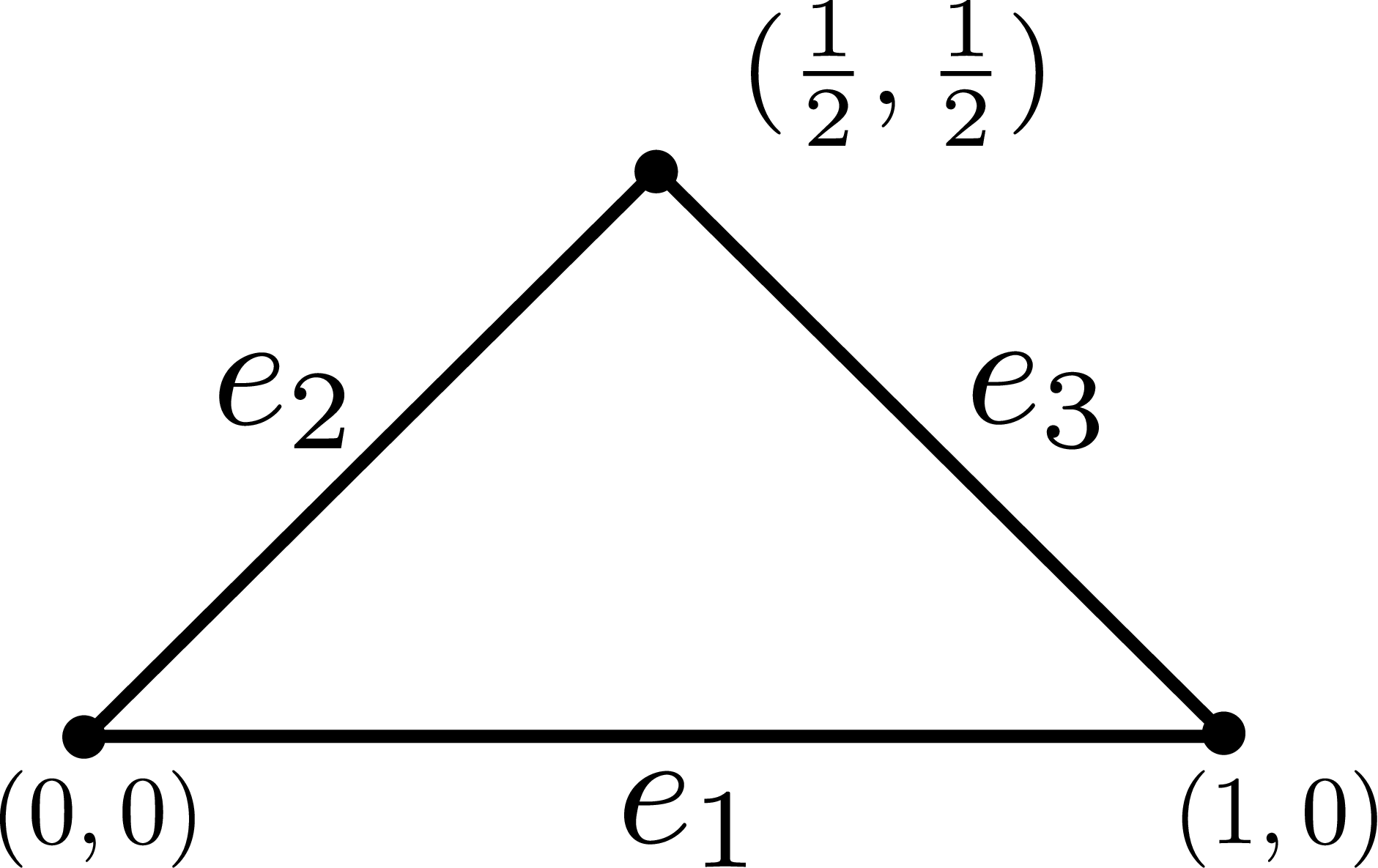}
   \includegraphics[width=3.6cm]{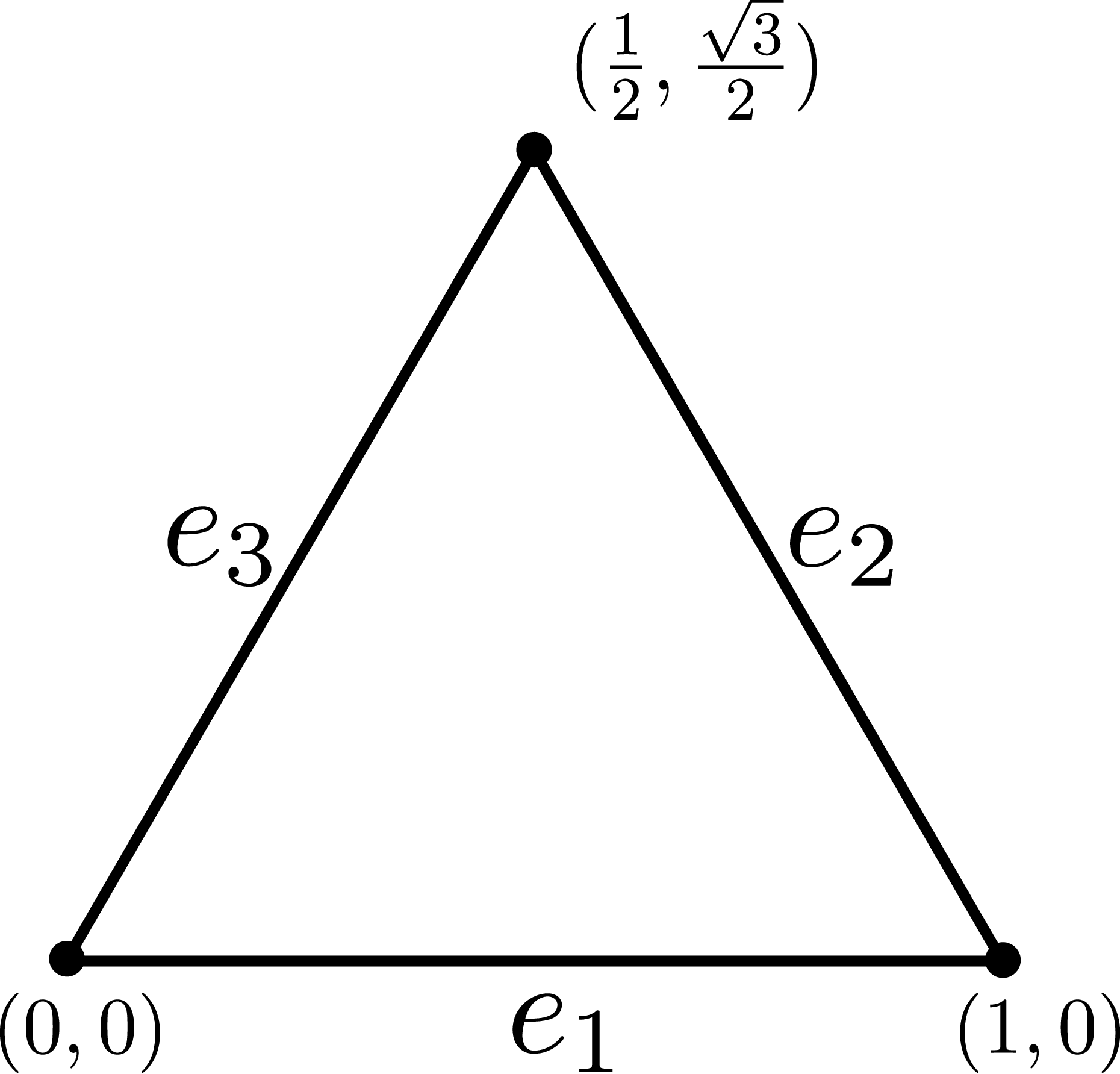}
   \includegraphics[width=3.6cm]{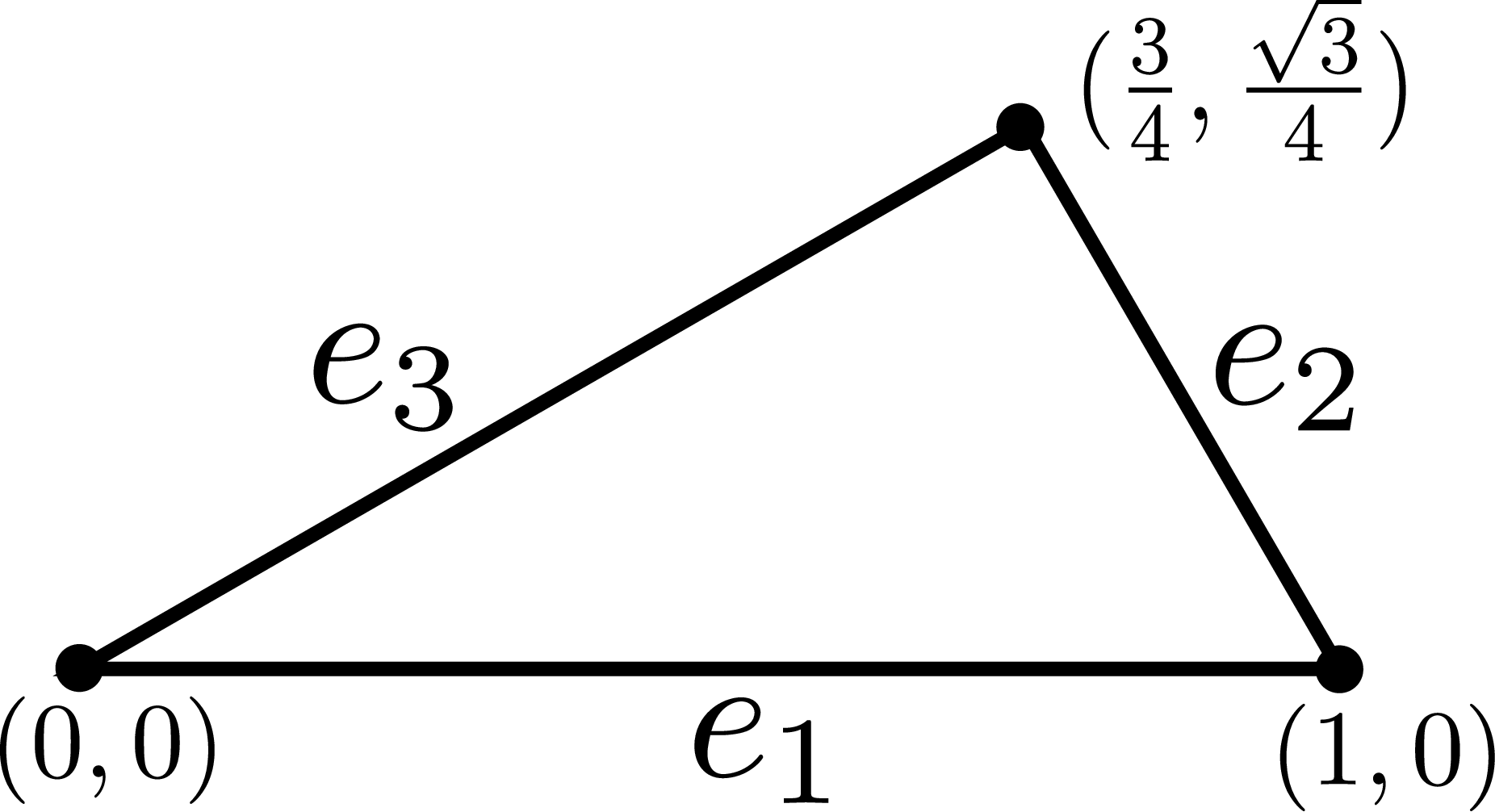}
   
   $K_1$ \qquad \quad \qquad \qquad \qquad \qquad $K_2$ \qquad \qquad \qquad \quad \qquad \qquad $K_3$
\end{center}
\caption{Three types of triangles
}
\label{fig02}
\end{figure}

 We consider three types of triangles (see  Fig. \ref{fig02}) and evaluate $C_{e_i}(K)$ by solving the corresponding Steklov eigenvalue problems using the linear conforming FEM. 
The results are shown in Tab. \ref{table:comparison_C_K}. 
It is observed that the bound in 
\eqref{eq:c_lliu} is not too {rough} and 
a direct estimation of $C_{e_i}(K)$ by solving the Steklov eigenvalue problem can obtain a sharper bound for the constant. 

\begin{table}[h!]
\begin{center}
\caption{Evaluation of $C_{e_i}(K)$ (mesh size $ h = 1/256$) \label{table:comparison_C_K}}
\begin{tabular}{c|ccc|ccc|ccc}
&
\multicolumn{3}{c|}{\rule[-2mm]{0mm}{6mm}{}
Approximation of $C_{e_i}(K)$} 
& \multicolumn{3}{c|}{Upper bound of $C_{e_i}(K)$} & \multicolumn{3}{c}{Upper bound in \eqref{eq:c_lliu}} \\ \cmidrule{2-10} 
& \rule[-2mm]{0mm}{5mm}{}$e_1$ & $e_2$ & $e_3$ & $e_1$ & $e_2$ & $e_3$ & $e_1$ & $e_2$ & $e_3$  \\ \midrule 
\rule[-1mm]{0mm}{5mm}{}
$K_1$ & 0.7071 & 0.5516 & 0.5516 & 0.7198 & 0.5571 & 0.5571 & 1.1481 & 0.9654 & 0.9654 \\ 
\rule[-1mm]{0mm}{5mm}{}
$K_2$ & 0.6361 & 0.6361 & 0.6361 & 0.6446 & 0.6446 & 0.6446 & 0.8723 & 0.8723 & 0.8723 \\ 
\rule[-1mm]{0mm}{5mm}{}
$K_3$ & 0.7700 & 0.4285 & 0.7071 & 0.7843 & 0.4320 & 0.7169 & 1.2337 & 0.8723 & 1.1480 \\ 
\end{tabular}
\end{center}
\end{table}

\section{ Conclusion}
In this paper, we propose a method to obtain the guaranteed lower bound of the Steklov eigenvalue by using the conforming FEM, where the hypercircle method plays an important role in obtaining the {\em a priori } error estimation. 
The proposed eigenvalue bounds have a degenerated convergence rate as $O(h)$, when the FEM approximations of the leading eigenvalues demonstrate the $O(h^2)$ convergence rate. 
{ Such a degenerated convergence rate of our propose method cannot be improved, because the involved projection error estimate has to handle the worst case when the solution to boundary value problem does not have the $H^2$ regularity. In future work, the authors will apply the Lehmann--Goerisch's theorem to obtain lower eigenvalue bounds with optimal convergence rates.}


\vspace{1cm}
\noindent \textbf{Funding:} The first author is supported by JST SPRING, Grant Number JPMJSP2121.
The second author has been supported by the National Natural Science Foundation of China (No.11426039,12061057,11571023).
The last author is supported by Japan Society for the Promotion of Science: Fund for the Promotion of Joint International Research (Fostering Joint International Research (A)) 20KK0306, Grant-in-Aid for Scientific Research (B) 20H01820, 21H00998.
This work also received support from the Research Institute for Mathematical Sciences, an International Joint Usage/Research Center located in Kyoto University.

\bibliographystyle{MyCMAM}
\bibliography{reference}   

\end{document}